\documentclass[10pt]{amsart}

\usepackage{amsmath}
\usepackage{amsthm}
\usepackage{amssymb}
\usepackage{amscd}
\usepackage{diagrams}

\theoremstyle{plain}
\newtheorem{dfn}[subsection]{Definition}
\newtheorem{thm}[subsection]{Theorem}
\newtheorem{prp}[subsection]{Proposition}
\newtheorem{cor}[subsection]{Corollary}
\newtheorem{lma}[subsection]{Lemma}
\newtheorem{sthm}[subsubsection]{Theorem}
\newtheorem{sprp}[subsubsection]{Proposition}
\newtheorem{slma}[subsubsection]{Lemma}

\theoremstyle{remark}
\newtheorem{rmk}[subsection]{Remark}
\newtheorem{exms}[subsection]{Examples}
\newtheorem{srmk}[subsubsection]{Remark}

\def\DD{\mathcal{D}}

\def\EE{\mathcal{E}}
\def\GG{\mathcal{G}}

\def\TT{\mathbb{T}}
\def\NN{\mathrm{N}}

\def\DD{\mathcal{D}}
\def\EE{\mathcal{E}}
\def\KK{\mathcal{K}}

\def\UU{\mathcal{U}}

\def\FF{\mathcal{F}}
\def\PP{\mathcal{P}}
\def\AA{\mathbb{A}}
\def\BB{\mathbb{B}}
\def\RR{\mathcal{R}}
\def\SS{\mathcal{S}}
\def\ZZ{\mathbb{Z}}
\def\WW{\mathrm{W}}

\def\nn{\underline{n}}

\def\Sing{\mathrm{Sing}}
\def\Coll{\mathrm{Coll}}
\def\Oper{\mathrm{Oper}}
\def\Cooper{\mathrm{Cooper}}

\def\Dop{{\Delta^\mathrm{op}}}

\def\FF{\mathcal{F}}

\def\Ch{\mathrm{Ch}}
\def\Sets{\mathrm{Sets}}

\def\op{\mathrm{op}}
\def\Top{\mathrm{Top}}

\def\Hom{\mathrm{Hom}}
\def\iHom{\underline{\Hom}}

\def\Sg{\Sigma}
\def\sg{\sigma}
\def\eps{\epsilon}
\def\al{\alpha}
\def\ito{\rightarrowtail}

\def\lra{\longrightarrow}
\def\llra{\overset{\sim}{\longleftrightarrow}}
\def\lrto{\leftrightarrows}
\def\dto{\rightrightarrows}
\def\inc{\hookrightarrow}
\def\eqv{\overset{\sim}{\longrightarrow}}

\def\cli{\varinjlim}

\def\Aut{\mathrm{Aut}}

\begin{document}
\title[The Boardman-Vogt resolution of operads]{The Boardman-Vogt resolution of operads in monoidal model categories}

\author{Clemens Berger and Ieke Moerdijk}

\date{23 September 2005 (revised version)}

\begin{abstract}We extend the $\WW$-construction of Boardman and Vogt to operads of an arbitrary monoidal model category with suitable interval, and show that it provides a cofibrant resolution for well-pointed $\Sg$-cofibrant operads. The standard simplicial resolution of Godement as well as the cobar-bar chain resolution are shown to be particular instances of this generalised $\WW$-construction.\end{abstract}

\maketitle

\section*{Introduction}

In \cite{BM}, sufficient conditions are given for the existence of a model structure on operads in an arbitrary (symmetric) monoidal model category. These conditions imply in particular that each operad may be resolved by a cofibrant operad. This general existence result leaves open the relation to various explicit resolutions of operads in literature, like the $\WW$-construction of Boardman and Vogt \cite{BV} for topological operads (special PROP's in their terminology), the cobar-bar resolution for chain operads \cite{GJ}, \cite{GK}, and the standard simplicial resolution of Godement \cite{G} arising from the free-forgetful adjunction between pointed collections and operads. 

The purpose of this article is to introduce a general, inductively defined $\WW$-construction, and to prove that it provides a functorial \emph{cofibrant} resolution for operads whose underlying collection is cofibrant and well-pointed. Our $\WW$-construction applies in an arbitrary monoidal model category which comes equipped with a suitable interval, and specialises to each of the above mentioned resolutions by an appropriate choice of ambient model category and interval. It is completely uniform, and has remarkable functorial and homotopical properties. 

The main idea of Boardman and Vogt was to enrich the free operad construction by assigning lengths to edges in trees. Surprisingly, all we need to carry out this topological idea in general, is an abstract interval with suitable algebraic and homotopical properties. We call the underlying algebraic object a \emph{segment} (i.e. an augmented associative monoid with absorbing element), and use the term \emph{interval} for segments which, in addition, induce cylinder-objects in Quillen's model-theoretic sense \cite{Q}. The real unit-interval equipped with its maximum operation is an example of such an interval for compactly generated spaces, while the standard representable $1$-simplex is an interval for simplicial sets. Moreover, the model category of chain complexes has an interval with segment structure because the normalised chain functor transfers the required structure from simplicial sets to chain complexes.

We show that for reduced chain operads, the $\WW$-construction is isomorphic to the \emph{cobar-bar resolution}, cf. Kontsevich and Soibelman \cite{KS}. The fact that such a resolution for chain operads exists goes back to Ginzburg and Kapranov \cite{GK}. Our treatment of the cobar-bar adjunction closely follows Getzler and Jones \cite{GJ}.

In an arbitrary monoidal model category, the adjunction between pointed collections and operads gives rise to a simplicial Godement resolution for operads. We show that this Godement resolution is also an instance of the $\WW$-construction with respect to a simplicial segment. This raises the problem of realising simplicial objects with algebraic structure in a monoidal model category. We address this problem in the Appendix where we introduce for this purpose the concept of a \emph{standard system of simplices}. Roughly speaking, this is a cosimplicial object which endows the monoidal model category with specific \emph{framings} (cf. \cite{DHK}, \cite{Hir}, \cite{Hov}) compatible with the monoidal structure, at least up to coherent homotopy. Compactly generated spaces, simplicial sets, symmetric spectra and chain complexes have such standard systems of simplices. Therefore, the Godement resolution of an operad realises in an appropriate way in all the cited cases.

The $1$-truncation of a standard system of simplices defines an interval with segment structure, so we get two a priori different resolutions for a given well-pointed $\Sg$-cofibrant operad: the Boardman-Vogt resolution and the realisation of the Godement resolution. We construct a comparison map from the former to the latter, and show that in the categories of compactly generated spaces, of simplicial sets and of symmetric spectra, this comparison map is an isomorphism; in the category of chain complexes, it is a weak equivalence of $\Sg$-cofibrant resolutions.\vspace{2ex}

The plan of this article is as follows:\vspace{1ex}

Section 1 describes the $\WW$-construction in an informal ``set-theoretical'' way so as to provide some intuitive background that should facilitate reading the more technical part of the article.\vspace{1ex}

Section 2 contains the basic definitions and properties of monoidal model categories. In particular, we fix our convention about cofibrant and $\Sg$-cofibrant operads, state the overall used telescope and patching lemmas of Reedy, and prove some lemmas in equivariant monoidal homotopy theory which we need later.\vspace{1ex}

Section 3 reviews the construction of the free operad generated by a pointed collection and fixes our convention about trees.\vspace{1ex}

Section 4 defines the $\WW$-construction as a sequential colimit of pushouts in the category of collections. The different layers of this filtered $\WW$-construction are determined by the number of internal edges of the indexing trees.\vspace{1ex}

Section 5 deduces the cofibrancy of the $\WW$-construction (for well-pointed $\Sg$-cofibrant operads) from a factorisation of the counit of the adjunction between pointed collections and operads.\vspace{1ex}

Section 6 describes the homotopical behaviour of the $\WW$-construction under maps of segments, maps of operads and under symmetric monoidal functors. \vspace{1ex}

Section 7 presents a \emph{relative $\WW$-construction} for operads under a fixed operad.\vspace{1ex}

Section 8 studies various instances of the $\WW$-construction: the classical construction of Boardman and Vogt for topological operads, the simplicial Godement resolution and the cobar-bar resolution for chain operads. \vspace{1ex}

The Appendix studies realisation functors for simplicial objects of a monoidal model category. We establish a correspondence between \emph{global} properties of the realisation functor and \emph{local} properties of the defining cosimplicial object.\vspace{1ex}

\thanks{{\emph{Acknowledgements:} We are grateful to Brooke Shipley, Paul Arne \O stv\ae r, Bruno Vallette, Rainer Vogt and an anonymous referee for several helpful remarks.}}

\section{Informal description of the $\WW$-construction}

By a \emph{segment} $H$ in a symmetric monoidal category $\EE$ with unit $I$, we mean an augmented associative monoid equipped with an absorbing element. We denote the binary operation by $\vee:H\otimes H\to H$, the neutral element by $0:I\to H$, the absorbing element by $1:I\to H$, and the augmentation by $\eps:H\to I$. Set-theoretically, the equations satisfied by the binary operation are the following : $x\vee(y\vee z)=(x\vee y)\vee z\,,\,0\vee x=x\vee 0\,,\,1\vee x=1=x\vee 1$. The reader may think of the real unit-interval equipped with the maximum operation.

Let $P$ be an operad in $\EE$. According to Boardman and Vogt \cite[Chapter III]{BV} (who consider topological operads only), a new operad $\WW(P)$ can be constructed as follows. Elements of $\WW(P)(n)$ are represented by planar trees, with exactly $n$ input edges; furthermore, each vertex of valence $k$ is labelled by an element of $P(k)$, and each internal edge is labelled by an element of $H$, which can be thought of as its ``length''. There are several types of identifications to be made, specifying when two such labelled trees define the same element in $\WW(P)(n)$.

One type of identification has to do with automorphisms of trees; we will be more explicit about this in Section $3$. Furthermore, there are two types of identifications relating to the operad composition of $P$: If an edge has length $0$, the tree is identified with the one obtained by contracting this edge and applying the corresponding circle-operation in $P$. If a vertex $v$ is labelled by the unit of the operad $P$, the tree is identified with the one obtained by deleting $v$, and by taking the maximum $t_1\vee t_2$ of the corresponding lengths of the edges incident to $v$ if these are internal. If $t_1$ or $t_2$ happens to be $0$, this identification is compatible with the previous one because $0$ is neutral. If one of the edges incident to $v$ is external, one simply deletes the length of the other edge.

In this way, one obtains a new operad $\WW(P)$, with operad composition defined by the usual grafting of trees, giving the new internal edges length $1$. This grafting operation is well defined on equivalence classes. In particular, the last identification above is compatible with grafting, precisely because $1$ is absorbing. Observe that the augmentation $H\to I$ induces a map of operads $\WW(P)\to P$, defined by forgetting lengths and applying composition in $P$. The free operad $\FF_*(P)$ generated by $P$ (viewed as a collection, pointed by the unit of $P$) can be identified with the suboperad of $\WW(P)$ given by the trees all of whose internal edges have length $1$. This defines a map $\FF_*(P)\to\WW(P).$ The composition $\FF_*(P)\to\WW(P)\to P$ is easily identified with the counit of the adjunction between pointed collections and operads. This factorisation has remarkable homotopical properties, whenever $\EE$ comes equipped with a compatible Quillen model structure and the segment $H$ induces cylinder-objects in the model-theoretic sense \cite[I.1]{Q}. Indeed, we prove in Theorem \ref{main} that under these hypotheses, if the collection underlying $P$ is cofibrant and well-pointed, the counit is factored into a \emph{cofibration} $\FF_*(P)\ito\WW(P)$ followed by a \emph{weak equivalence} $\WW(P)\eqv P$. Since under the same hypotheses $\FF_*(P)$ is a cofibrant operad, the $\WW$-construction thus provides a \emph{cofibrant resolution} for $P$.\vspace{2ex}

\section{Monoidal model categories}

 \subsection{Basic definitions}\label{ppl}We will work in an arbitrary \emph{monoidal model category} $\EE$. This means that $\EE$ is a closed model category equipped with a compatible monoidal structure, which we always assume to be \emph{symmetric} and \emph{closed}, cf. \cite{Hov}, \cite{SS}. We write $\otimes$ for the symmetric monoidal structure, $I$ for the unit, and $Y^X$ for the internal hom of two objects $X,Y$ of $\EE$. The compatibility axiom is the so-called \emph{pushout-product axiom}, which states that for any two cofibrations $A\ito B$ and $X\ito Y$, the induced map $(A\otimes Y)\cup_{A\otimes X}(B\otimes X)\to B\otimes Y$ is again a cofibration, which is trivial if $A\ito B$ or $X\ito Y$ is.

Given $n$ cofibrations $X^i_0\ito X^i_1,\,1\leq i\leq n$, the tensors $X^1_{\eps_1}\otimes\cdots\otimes X^n_{\eps_n}$, with $\eps_i\in\{0,1\}$ and \emph{not} all $\eps_i$ equal to $1$, fit into a diagram whose colimit maps to $X^1_1\otimes\cdots\otimes X_1^n$. A repeated application of the pushout-product axiom shows that this map is a cofibration, which is trivial if one of $X^i_0\ito X^i_1$ is. We refer to this more general property for $n$ cofibrations as the \emph{pushout-product lemma}.

Throughout, we assume that the unit $I$ of $\EE$ is \emph{cofibrant}. If the latter is not the case, one should add a further axiom, the so-called \emph{unit axiom}, which guarantees that cofibrant resolutions of the unit are compatible with the monoidal structure. In our context, we need the cofibrancy of the unit because otherwise the concept of an interval, see Definition \ref{interval}, would be much harder to handle.

An adjunction between model categories is a \emph{Quillen adjunction} if the left adjoint preserves cofibrations and the right adjoint preserves fibrations. The left adjoint of a Quillen adjunction is often called a \emph{left Quillen functor}, and the right adjoint a \emph{right Quillen functor}. By adjointness, any left Quillen functor preserves trivial cofibrations, and any right Quillen functor preserves trivial fibrations, cf. \cite[8.5.3]{Hir}. K. Brown's lemma \cite[1.1.12]{Hov} implies then that any left Quillen functor also preserves weak equivalences between cofibrant objects, while any right Quillen functor also preserves weak equivalences between fibrant objects, cf. \cite[8.5.7]{Hir}.

\subsection{Symmetric monoidal and h-monoidal functors}\label{functor}A functor  between symmetric monoidal categories $\phi:(\EE,I,\otimes,\tau)\to(\EE',I',\otimes',\tau')$ is \emph{symmetric monoidal} if $\phi$ comes equipped with a unit map $\phi_0:I'\to\phi(I)$ and with binatural maps $\phi_{XY}:\phi(X)\otimes'\phi(Y)\to\phi(X\otimes Y)$ satisfying well known associativity, symmetry and unit constraints, cf. \cite[Chapter 20]{MMSS}, \cite{KMac}. The functor $\phi$ is \emph{strong symmetric monoidal} if $\phi_0$ and all $\phi_{XY}$ are isomorphisms. Symmetric monoidal functors preserve algebraic (though not coalgebraic) structures, in particular they preserve monoids, commutative monoids, operads, and algebras over operads. 

A functor $\phi:\EE\to\EE'$ between monoidal model categories will be called \emph{h-monoidal} if it is a symmetric monoidal functor such that $\phi_0$ is an isomorphism, and $\phi_{XY}$ are \emph{weak equivalences} whenever $X$ and $Y$ are \emph{cofibrant}.

\subsection{Reedy's patching and telescope lemmas}

In proving that a certain map is a weak equivalence, we shall frequently use a device which appeared (in this generality) for the first time in Chris Reedy's PhD thesis, cf. \cite{Hir}, \cite{Hov}: assume that our map is a colimit $\cli\phi$ of a natural transformation $\phi:F\overset{\cdot}{\to}G$ of $\DD$-diagrams $F,G:\DD\to\EE$ in a model category $\EE$. Such a colimit is a weak equivalence whenever the following three conditions are satisfied:\begin{enumerate}\item $\DD$ is equipped with the structure of a direct Reedy category;\item For all objects $d$ of $\DD$, $\phi(d):F(d)\to G(d)$ is a weak equivalence;\item $F$ and $G$ are cofibrant with respect to the Reedy model structure on $\EE^\DD$.\end{enumerate}

We recall from \cite[12.8]{DHK} that a Reedy category $\DD$ is \emph{direct} if the inverse subcategory of $\DD$ (which contains all objects, but only those morphisms that lower degrees) is a coproduct of categories with terminal object. For any direct Reedy category $\DD$, the colimit functor $\cli:\EE^\DD\to\EE$ is a left Quillen functor with respect to the Reedy model structure on $\EE^\DD$ since, under the above condition on $\DD$, the right adjoint diagonal $\EE\to\EE^\DD$ preserves fibrations and trivial fibrations. By (2) and (3), the natural transformation $\phi:F\overset{\cdot}{\to}G$ is a weak equivalence between cofibrant objects in $\EE^\DD$. K. Brown's lemma implies then that $\cli\phi:\cli F\to\cli G$ is a weak equivalence between cofibrant objects in $\EE$.

\emph{Reedy's patching lemma} is the above statement for the direct Reedy category $\DD_1=\cdot\leftarrow\cdot\rightarrow\cdot$, where we assume that the left arrow lowers degree, while the right arrow raises degree. A $\DD_1$-diagram is Reedy-cofibrant if the three objects are cofibrant, and the arrow on the right-hand side is a cofibration. \emph{Reedy's telescope lemma} is the above statement for the direct Reedy category $\DD_2=(\mathbb{N},<)$. A $\DD_2$-diagram is Reedy-cofibrant if all objects are cofibrant and all arrows are cofibrations.

\subsection{Cofibrant and $\Sg$-cofibrant operads}\label{cofoperad}

Let $\EE$ be a \emph{cofibrantly generated} mon-oidal model category, cf. \cite{Hir}, \cite{Hov}, \cite{SS}. For any discrete group $\Gamma$, the monoidal model structure on $\EE$ can then be transferred to a monoidal model structure on the category $\EE^\Gamma$ of objects in $\EE$ equipped with a right $\Gamma$-action. The \emph{monoidal structure} is determined by letting $\Gamma$ act diagonally on tensor products, and by conjugation on internal hom's. The \emph{model structure} is determined by the property that a map in $\EE^\Gamma$ is a fibration (resp. weak equivalence) if and only if the underlying map in $\EE$ is a fibration (resp. weak equivalence). 

We will refer to a cofibration in $\EE^\Gamma$ as a $\Gamma$-cofibration in $\EE$, and to a cofibrant object in $\EE^\Gamma$ as a $\Gamma$-cofibrant object in $\EE$. For the permutation groups $\Sg_n$, these model structures on $\EE^{\Sg_n}$ together define a model structure on the \emph{category of collections},$$\Coll(\EE)=\prod_{n\geq 0}\EE^{\Sg_n}.$$

We write $\Oper(\EE)$ for the \emph{category of operads} in $\EE$, and call a map $Y\to X$ of operads a fibration (resp. weak equivalence) if the underlying map of collections is a fibration (resp. weak equivalence); i.e. if for each $n$, the map $Y(n)\to X(n)$ is a fibration (resp. weak equivalence) in $\EE$. Under general conditions, this defines a model structure on the category $\Oper(\EE)$, see \cite{BM,Sp}. 

\emph{Independently of this}, we call a map $P\to Q$ of operads a \emph{cofibration} if it has the left lifting property with respect to the trivial fibrations $Y\to X$ of operads just defined, and a \emph{$\Sg$-cofibration} if the underlying map of collections is a cofibration in the model category of collections. In particular, an operad $P$ is called \emph{cofibrant} (resp. \emph{$\Sigma$-cofibrant)} if the unique map from the initial operad to $P$ is a cofibration (resp. $\Sg$-cofibration) in the sense just described. Observe that a symmetric monoidal left Quillen functor preserves $\Sigma$-cofibrations of operads, but not necessarily cofibrations of operads, unless the right adjoint functor is also symmetric monoidal.

\subsection{Equivariant monoidal homotopy theory}

The homotopy theory of symmetric operads in a monoidal model category $\EE$ relies on some basic equivariant homotopy theory. In this section we collect those results which are needed later. \begin{slma}\label{eq1}For a cofibrantly generated model category $\EE$, induction and restriction along a group homomorphism $\Gamma_1\to\Gamma_2$ form a Quillen adjunction $\EE^{\Gamma_1}\leftrightarrows\EE^{\Gamma_2}$. Restriction along an inclusion $\Gamma_1\inc\Gamma_2$ takes $\Gamma_2$-cofibrations to $\Gamma_1$-cofibrations.\end{slma}

\begin{proof}Restriction takes fibrations (resp. weak equivalences) in $\EE^{\Gamma_2}$ to fibrations (resp. weak equivalences) in $\EE^{\Gamma_1}$; therefore, induction takes $\Gamma_1$-cofibrations to $\Gamma_2$-cofibrations, which proves the first part. For the second part, observe that restriction preserves all colimits (being also a left adjoint), and that the generating $\Gamma_2$-cofibrations split after restriction into a coproduct of generating $\Gamma_1$-cofibrations.\end{proof}

For the rest of this section we will use the following terminology: any map in $\EE^\Gamma$ whose underlying map is a cofibration in $\EE$ will be called a \emph{$\Gamma$-equivariant cofibration}.

\begin{slma}\label{eq2}Let $A\ito B$ and $X\ito Y$ be $\Gamma$-equivariant cofibrations. If one of them is a $\Gamma$-cofibration, the pushout-product map $(A\otimes Y)\cup_{A\otimes X}(B\otimes X)\to B\otimes Y$ is a $\Gamma$-cofibration. Moreover, the latter is trivial if $A\ito B$ or $X\ito Y$ is.\end{slma}

\begin{proof}Let $Z\to W$ be a $\Gamma$-equivariant map. By exponential transpose, the following two commutative diagrams in $\EE^\Gamma$ correspond to each other:\begin{diagram}[small,silent](A\otimes Y)\cup_{A\otimes X}(B\otimes X)&\rTo&Z&\quad\quad&A&\rTo&Z^Y\\\dTo&&\dTo&&\dTo&&\dTo\\B\otimes Y&\rTo&W&&B&\rTo&W^Y\times_{W^X}Z^X\end{diagram} In particular, the pushout-product map is a $\Gamma$-cofibration whenever, for any trivial fibration $Z\to W$ in $\EE^\Gamma,$ the right hand square admits a diagonal filler. If $A\to B$ is a $\Gamma$-cofibration, such a filler exists whenever $Z^Y\to W^Y\times_{W^X}Z^X$ is a trivial fibration in $\EE^\Gamma$ or, equivalently, in $\EE$. The latter holds by the pushout-product axiom, since $X\to Y$ is supposed to be a $\Gamma$-equivariant cofibration. A symmetric argument gives the result for a $\Gamma$-equivariant cofibration $A\to B$ and a $\Gamma$-cofibration $X\to Y$. The proof of the second statement follows the same pattern.\end{proof}

Lemma \ref{eq2} is the equivariant form of the pushout-product axiom. Like for the latter, there is a corresponding \emph{equivariant pushout-product lemma} for finite families of $\Gamma$-equivariant cofibrations.

\begin{slma}\label{eq3}Let $1\to\Gamma_1\to\Gamma\to\Gamma_2\to 1$ be a short exact sequence of groups. Let $A\ito B$ be a $\Gamma_2$-cofibration, and $X\ito Y$ be a $\Gamma$-equivariant $\Gamma_1$-cofibration. Then the pushout-product map $(A\otimes Y)\cup_{A\otimes X}(B\otimes X)\to B\otimes Y$ is a $\Gamma$-cofibration. Moreover, the latter is trivial if $A\ito B$ or $X\ito Y$ is.\end{slma}
\begin{proof}We consider an arbitrary trivial fibration $Z\to W$ in $\EE^\Gamma$. Since $\Gamma_1$ acts trivially on $A\ito B$, any map $A\to Z^Y$ (resp. $B\to W^Y\times_{W^X}Z^X$) factors through the fixed point object, so we have the following commutative diagram:\begin{diagram}[small,silent]A&\rTo&(Z^Y)^{\Gamma_1}&\rTo&Z^Y\\\dTo&&\dTo&&\dTo\\B&\rTo&(W^Y\times_{W^X}Z^X)^{\Gamma_1}&\rTo&W^Y\times_{W^X}Z^X\end{diagram}for which the existence of a $\Gamma$-equivariant diagonal filler has to be shown, like in the proof of Lemma \ref{eq2}. For this, it suffices to construct a $\Gamma_2$-equivariant diagonal filler for the left hand square. Since $A\ito B$ is a $\Gamma_2$-cofibration, it remains to be shown that the middle vertical map is a trivial fibration in $\EE^{\Gamma_2}$ or, equivalently, in $\EE$. But the latter follows, again by exponential transpose, from Lemma \ref{eq2}, using that $X\ito Y$ is a $\Gamma_1$-cofibration as well as the fact that the fixed point functor $(-)^{\Gamma_1}$ has a left adjoint. The proof of the second statement follows the same pattern.\end{proof}

Lemma \ref{eq3} recovers Lemma 5.10 of \cite{BM} which treats the special case of a split exact sequence $1\to G\to G\rtimes\Sg\to\Sg\to 1$. The diligent reader has certainly observed that in the statement of the latter a $G$ is missing. 

The following two lemmas are needed in the proof of Propositions \ref{cofibrant1}, \ref{cofibrant2} and \ref{relative}. In their statement, we use the union symbol to indicate specific colimits in $\EE$. More precisely, $n$-fold tensor products of sequences $A_i^1\to A_i^2\to A_i^3,\,1\leq i\leq n$, fit into a subdivided $n$-cubical diagram, with vertices $A_1^{i_1}\otimes A_2^{i_2}\otimes\cdots\otimes A^{i_n}_n$ where $i_k\in\{1,2,3\}$ and $1\leq k\leq n$. Any set $\mathcal{J}$ of such vertices defines a cone consisting of the objects of $\mathcal{J}$ and all arrows and vertices \emph{above} objects of $\mathcal{J}$. For $\mathcal{J}=\{S_1,\dots,S_k\}$, we write $S_1\cup\cdots\cup S_k$ to indicate the colimit of this cone; the colimit comes equipped with a canonical map to $A_1^3\otimes\cdots\otimes A_n^3$.

\begin{slma}\label{eq4}Let $A_i\ito B_i\ito C_i$ be $\Gamma$-equivariant cofibrations for $1\leq i\leq n$. Then, the induced map$$\bigotimes_{i=1}^nB_i\cup\bigcup_{i=1}^n(C_1\otimes\cdots\otimes C_{i-1}\otimes A_i\otimes C_{i+1}\otimes\cdots\otimes C_n)\to C_1\otimes\cdots\otimes C_n$$ is a $\Gamma$-equivariant cofibration. Moreover, the latter is trivial (resp. a $\Gamma$-cofibration) if all $B_i\ito C_i$ are trivial (resp. if for some $i$, $A_i\ito B_i\ito C_i$ are $\Gamma$-cofibrations).\end{slma}

\begin{proof}For $n=1$, the lemma is true. Assume the lemma is proved for $n-1$ sequences $A_i\ito B_i\ito C_i$. For shortness, we shall suppress tensors from notation. The induced map $\gamma_n$ for $n$ sequences is then the composite of $\alpha_n$ and $\beta_n$, as indicated in the following commutative triangle:\begin{diagram}[silent]B_1\cdots B_n\cup\bigcup_{i=1}^nC_1\cdot\cdot C_{i-1}A_iC_{i+1}\cdot\cdot C_n&&\\\dTo^{\al_n}&\rdTo^{\gamma_n}&\\B_1\cdot\cdot B_{n-1}C_n\cup\bigcup_{i=1}^nC_1\cdot\cdot C_{i-1}A_iC_{i+1}\cdot\cdot C_n&\rTo^{\beta_n}&C_1\cdots C_n\end{diagram}The map $\al_n$ is part of the pushout\begin{diagram}[silent]B_1\cdots B_n\cup\bigcup_{i=1}^{n-1}B_1\cdot B_{i-1}A_iB_{i+1}\cdot B_{n-1}C_n&\rTo&B_1\cdots B_n\cup\bigcup_{i=1}^nC_1\cdot C_{i-1}A_iC_{i+1}\cdot C_n\\\dTo&&\dTo_{\al_n}\\B_1\cdot B_{n-1}C_n&\rTo&B_1\cdot B_{n-1}C_n\cup\bigcup_{i=1}^nC_1\cdot C_{i-1}A_iC_{i+1}\cdot C_n\end{diagram}the left vertical map of which is a $\Gamma$-equivariant cofibration by the pushout-product lemma (cf. \ref{ppl}) applied to $A_i\ito B_i,1\leq i\leq n-1,$ and $B_n\ito C_n$. Therefore, $\al_n$ is a $\Gamma$-equivariant cofibration. The map $\beta_n$ is obtained as the pushout-product map of $\gamma_{n-1}$ and $A_n\ito C_n$, and is therefore also a $\Gamma$-equivariant cofibration.

If all $B_i\ito C_i$ are trivial, then so is $\gamma_n=\beta_n\al_n$ by the same argument as above. If there is some $i$ for which $A_i\ito B_i\ito C_i$ are $\Gamma$-cofibrations, then, using an equivariant form of the pushout-product lemma based on Lemma \ref{eq2}, the preceding proof actually shows that $\gamma_n=\beta_n\al_n$ is a $\Gamma$-cofibration.\end{proof}

\begin{slma}\label{eq5}Let $1\to\Gamma_1\to\Gamma\to\Gamma_2\to 1$ be a short exact sequence of groups. Let $A\ito B\ito C$ be $\Gamma_2$-cofibrations, and let $X\ito Y\ito Z$ be $\Gamma$-equivariant $\Gamma_1$-cofibrations. Then the induced maps$$(A\otimes Z)\cup(B\otimes Y)\cup(C\otimes X)\to(B\otimes Z)\cup(C\otimes X)\to C\otimes Z$$are $\Gamma$-cofibrations. Moreover, the latter are trivial if $B\ito C$ and $Y\ito Z$ are.\end{slma}

\begin{proof}We argue like in the proof of Lemma \ref{eq4} for $n=2$, using Lemma \ref{eq3} instead of Lemma \ref{eq2}.\end{proof}

\section{Free operads and trees}

The Boardman-Vogt resolution is a parametrised form of the free operad construction, and like the latter is defined in terms of trees. Therefore,  it is necessary to be more explicit about these first. As in \cite{BM}, we write $\TT$ for the groupoid of planar trees and non-planar isomorphisms. By a planar tree, we mean a finite planar tree with one output edge on the bottom, and input edges on the top. For example,
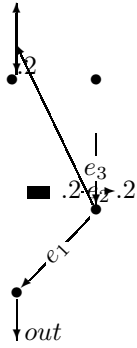
\begin{diagram}[height=0.55cm,abut]&&&&&\\\rdTo(0.9,2)\dTo\hspace{0.4em}\ldTo(1.2,2)&&\\\bullet\,\,&\bullet&&\\&\rdTo(1.2,2)~{e_2}\dTo~{e_3}\ldTo(1.2,2)&\\&\bullet&&&\\\rdTo(1,2)\ldTo(1,2)~{e_1}&&&\\\bullet&&&\\\dTo_{out}&&&&\end{diagram}is a tree with $5$ input edges. These input edges and the edge from the root are the \emph{external} edges of the tree. The other edges ($e_1$, $e_2$ and $e_3$) are called \emph{internal}. There are four vertices here. The planar structure of the tree defines a linear order on the set of edges coming into each vertex. The cardinality of this set is called the valence of the vertex (here, the root has valence $2$, two vertices have valence three, and one vertex has valence $0$).

For a tree $T$, we write $in(T)$ for the set of its input edges, and $\lambda(T)$ for the set of orderings $\{1,\dots,n\}\overset{\sim}{\lra}in(T)$ where $n$ is the cardinality of $in(T)$. Notice that $\Sg_n$ acts on $\lambda(T)$ from the right in an evident way. Observe that the input edges of $T$ are linearly ordered by the planar structure of $T$, so that we can identify $\lambda(T)$ with the permutation group $\Sg_n$. Let $\TT(n)$ be the subgroupoid of $\TT$ consisting of trees with exactly $n$ input edges. Then $\lambda$ defines a functor $\lambda_n:\TT(n)\lra\Sets^{\Sg_n}.$ The colimit $\mathcal{T}(n)=\varinjlim\lambda_n$ is the $\Sg_n$-set of planar trees with $n$ ordered input edges. 

The collection $(\mathcal{T}(n))_{n\geq 0}$ underlies an operad, the well known \emph{operad of planar trees}. The operad composition is defined in terms of grafting of trees; the unit for this operation is the tree $|$ in $\TT(1)$ without vertices, in which the input edge is also the output edge. For each $n\geq 0$, we write $t_n$ for the tree with one vertex and $n$ inputs. If $T\in\TT$ is any tree other than the unit tree $|$, and $k$ is the valence of the root vertex of $T$, then $T$ is obtained by grafting trees $T_1,\dots,T_k$ onto $t_k$, formally$$T=t_k(T_1,\dots,T_k).$$This allows induction on trees. For instance, the automorphism group $\Aut(T)$ can be described as follows: up to isomorphism, $T$ is obtained as\begin{gather}\label{semi1}T=t_k(T_1^1,\dots,T_{k_1}^1,T_1^2,\dots,T_{k_2}^2,\dots,T_1^l,\dots,T_{k_l}^l)\end{gather}where $T_1^i,\dots,T_{k_i}^i$ are copies of one and the same tree $T^i$, and $T^i$ is not isomorphic to $T^j$ if $i\not=j$. Then an automorphism of $T$ maps each copy $T^i_r$ to another copy $T^i_{r'}$ of $T^i$ via an automorphism of $T^i$. Thus the automorphism group of $T$ is a semi-direct product\begin{gather}\label{semi2}\Aut(T)\cong(\Aut(T^1)^{k_1}\times\cdots\times\Aut(T^l)^{k_l})\rtimes(\Sg_{k_1}\times\cdots\times\Sg_{k_l})=\Gamma_T\rtimes\Sg_T\end{gather}where $\Sg_{k_i}$ acts on the product $\Aut(T^i)^{k_i}$ by permuting the factors, cf. \cite[5.8-9]{BM}.\vspace{1ex}

Let $\EE$ be  a monoidal model category as above and let $K=(K(n))_{n\geq 0}$ be a collection in $\EE$. We begin by giving an explicit description of the free operad on $K$, in a way which closely follows \cite{BM}. The collection $K$ gives rise to a functor$$\underline{K}:\TT^\op\to\EE$$defined inductively as follows. On objects, we set $\underline{K}(|)=I$ and$$\underline{K}(t_n(T_1,\dots,T_n))=K(n)\otimes\underline{K}(T_1)\otimes\cdots\otimes\underline{K}(T_n).$$To define $\underline{K}$ on arrows, note that if $\phi:T\to T'$ is an isomorphism, we can write $T=t_n(T_1,\dots,T_n)$ and $T'=t_n(T'_1,\dots,T'_n)$, and $\phi$ is given by a permutation $\sg\in\Sg_n$ together with isomorphisms $\phi_i:T_i\to T'_{\sg(i)}\quad(i=1,\dots,n)$. The map $\underline{K}(\phi):\underline{K}(T')\to\underline{K}(T)$ is then defined as the composition of the canonical isomorphism $$\sg^*\otimes\sg:K(n)\otimes\underline{K}(T'_1)\otimes\cdots\otimes\underline{K}(T'_n)\overset{\sim}{\lra}K(n)\otimes\underline{K}(T'_{\sg(1)})\otimes\cdots\otimes\underline{K}(T'_{\sg(n)})$$given by the right action of $\sg$ on $K(n)$ and the symmetry of $\EE$, and the map $id_{\underline{K}(t_n)}\otimes\underline{K}(\phi_1)\otimes\cdots\otimes\underline{K}(\phi_n)$. Then the free operad $\FF(K)$ can be computed as$$\FF(K)(n)=\coprod_{[T],T\in\TT(n)}\underline{K}(T)\otimes_{\Aut(T)}I[\Sg_n].$$Here, $[T]$ ranges over isomorphism classes of trees in $\TT(n)$, and $I[\Sg_n]=\coprod_{\sg\in\Sg_n}I$ with its canonical $\Sg_n$-action. Then $\Aut(T)$ acts on $I[\Sg_n]$ by the identification of $\lambda(T)$ with $\Sg_n$ pointed out above.

If $K$ is a cofibrant collection, then so is $\FF(K)$; in other words, $\FF(K)$ is a $\Sg$-cofibrant operad. In fact, for each tree $T$, the object $\underline{K}(T)$ is $\Aut(T)$-cofibrant. To see this, we use induction on trees and notation (\ref{semi1}) and (\ref{semi2}). By induction, each $\underline{K}(T^i_j)$ is $\Aut(T^i)$-cofibrant. It follows that $\underline{K}(T^1)^{\otimes k_1}\otimes\cdots\otimes\underline{K}(T^l)^{\otimes k_l}$ is $\Gamma_T$-cofibrant and has an action by $\Gamma_T\rtimes\Sg_T$. Also, $K(n)$ is $\Sg_n$-cofibrant by assumption, hence $\Sg_T$-cofibrant by Lemma \ref{eq1}. Lemma \ref{eq3}, applied to the exact sequence $1\to\Gamma_T\to\Gamma_T\rtimes\Sg_T\to\Sg_T\to 1$, then completes the argument.

A \emph{pointed} collection is a collection $K$ equipped with a base-point $I\to K(1)$. We write $\Coll_*(\EE)$ for the category of pointed collections and pointed maps between them. Since each operad has a unit, there is a forgetful functor$$\Oper(\EE)\to\Coll_*(\EE).$$This functor has a left adjoint,$$\FF_*:\Coll_*(\EE)\to\Oper(\EE),$$which associates to each pointed collection $K$ an operad $\FF_*(K)$. The feature that distinguishes it from $\FF(K)$ is that the base-point of $K$ is turned into the unit of the operad $\FF_*(K)$, whereas the unit of $\FF(K)$ is added ``freely''. The construction of $\FF_*(K)$ is more relevant to our present purposes than that of $\FF(K)$. Therefore we consider it in some detail. Formally, $\FF_*(K)$ is defined by the following pushout of operads\begin{diagram}[small,silent]\FF(\underline{I})&\rTo&\FF(K)\\\dTo&&\dTo\\\underline{I}&\rTo&\FF_*(K)\end{diagram}where $\underline{I}$ denotes the initial operad in $\EE$, i.e. $\underline{I}(1)=I$ and $\underline{I}(n)=0$ for $n\not=1$. Thus, in the terminology of \cite{BM}, $\FF_*(K)$ is the cellular extension of the initial operad with respect to the base-point $I\to K(1)$ of $K$, provided that this base-point is a cofibration. Pointed collections $K$ such that $I\to K(1)$ is a cofibration, will henceforth be called \emph{well-pointed}. In other words, for a pointed collection $K$, the canonical map of collections $\underline{I}\to K$ is a cofibration of collections if and only if $K$ is well-pointed and cofibrant.

\begin{lma}\label{cofibrant}For any well-pointed cofibrant collection $K$ of a cofibrantly generated monoidal model category $\EE$, the operad $\FF_*(K)$ is cofibrant.\end{lma}

\begin{proof}Since $\underline{I}\to\KK$ is a cofibration of collections, its image $\FF(\underline{I})\to\FF(\KK)$ under the free functor is a cofibration of operads (cf. \ref{cofoperad}). Therefore, the operadic pushout $\underline{I}\to\FF_*(\KK)$ is a cofibration of operads, and hence $\FF_*(\KK)$ is a cofibrant operad.\end{proof}  

For later purposes, it will be useful to give an explicit description of $\FF_*(K)$. We assume that $K$ is well-pointed and cofibrant (although this assumption is not necessary for the construction itself). For a tree $T\in\TT$, a vertex $v\in T$ is called \emph{unary} if it has valence $1$. For a set $c$ of unary vertices in $T$, we shall define an object $\underline{K}_c(T)$. Intuitively speaking, $\underline{K}(T)$ consists of all assignments $p$ of elements $p(v)\in K(|v|)$ to vertices $v\in T$ of valence $|v|$; and $\underline{K}_c(T)$ is the subobject defined by the requirement that $p(v)$ is the unit in $K(1)$ whenever $v\in c$. Formally, $\underline{K}_c(T)$ is defined by induction on $T$. We set $\underline{K}_c(|)=I$ and, if $T=t_n(T_1,\dots,T_n)$,$$\underline{K}_c(T)=\begin{cases}I\otimes\underline{K}_c(T_1)\otimes\cdots\otimes\underline{K}_c(T_n)\text{ if the root vertex of }T\text{ belongs to }c;\\ K(n)\otimes\underline{K}_c(T_1)\otimes\cdots\otimes\underline{K}_c(T_n)\text{ otherwise}.\end{cases}$$For non-empty sets $c,d$ of unary vertices in $T$ such that $c\subseteq d$, the unit $u:I\to K(1)$ induces a cofibration $\underline{K}_d(T)\to\underline{K}_c(T)$ by the pushout-product axiom. Notice that $\underline{K}_\emptyset(T)=\underline{K}(T)$. Let us write\begin{gather}\label{minus1}\underline{K}^-(T)=\bigcup_{\emptyset\not=c}\underline{K}_c(T)\end{gather}where $c$ ranges over all non-empty sets of unary vertices in $T$, and the union is interpreted as the colimit over all cofibrations $\underline{K}_d(T)\to\underline{K}_c(T)$ for $c\subseteq d$. An equivariant form of the pushout-product lemma, based on Lemma \ref{eq2}, together with an induction on trees shows then that the induced map \begin{gather}\label{minus2}\underline{K}^-(T)\to\underline{K}(T)\end{gather} is an $\Aut(T)$-cofibration for every well-pointed cofibrant collection $K$.

For any non-empty set $c$ of unary vertices in $T$, let us write $T/c$ for the tree obtained from $T$ by deleting each vertex in $c$ and joining the corresponding edges. Notice that $T/c$ has strictly less internal edges than $T$. Notice also that there is an evident (iso)morphism\begin{gather}\label{minus3}\underline{K}_c(T)\overset{\sim}{\lra}\underline{K}(T/c)\end{gather}We are now ready to construct the operad $\FF_*(K)$. Each $\FF_*(K)(n)$ is constructed as the colimit of a sequence of $\Sg_n$-cofibrations$$\FF_*(K)(n,0)\ito\FF_*(K)(n,1)\ito\cdots$$Intuitively speaking, $\FF_*(K)(n,k)$ is to be the part of $\FF_*(K)(n)$ constructed by trees with at most $k$ unary vertices. Let $\TT(n,k)$ be the subgroupoid of $\TT(n)$ consisting of trees with $n$ input edges and $k$ unary vertices. Let$$\FF_*(K)(n,0)=\coprod_{[T],T\in\TT(n,0)}\underline{K}(T)\otimes_{\Aut(T)}I[\Sg_n]$$Here, $[T]$ ranges over isomorphism classes of trees in $\TT(n,0)$. So, $\FF_*(K)(n,0)$ is a summand of $\FF(K)(n)$. Suppose $\FF_*(n,k-1)$ has been constructed and form the pushout\begin{diagram}[small,silent]\coprod_{[T],T\in\TT(n,k)}\underline{K}^-(T)\otimes_{\Aut(T)}I[\Sg_n]&\rTo&\FF_*(K)(n,k-1)\\\dTo&&\dTo\\\coprod_{[T],T\in\TT(n,k)}\underline{K}(T)\otimes_{\Aut(T)}I[\Sg_n]&\rTo&\FF_*(K)(n,k)\end{diagram}where the left vertical map is given by (\ref{minus2}) and the top horizontal map by (\ref{minus3}). Let$$\FF_*(K)(n)=\varinjlim_k\FF_*(K)(n,k).$$Notice that, since each of the maps (\ref{minus2}) is an $\Aut(T)$-cofibration, the vertical maps above are $\Sg_n$-cofibrations. In particular, $\FF_*(K)$ is $\Sg$-cofibrant.

\section{The $\WW$-construction}

In this section we give an inductive definition of the Boardman-Vogt resolution $\WW(H,P)$ of an operad $P$ (cf. Section $1$), in the general context of a cofibrantly generated monoidal model category $\EE$ (cf. Section $2$) with interval $H$. Recall that the unit $I$ of $\EE$ is cofibrant by assumption.

\begin{dfn}\label{interval}A \emph{segment} in $\EE$ is a factorisation $I\sqcup I\overset{(0,1)}{\ito}H\overset{\eps}{\lra}I$ of the codiagonal, together with an associative operation $\vee:H\otimes H\to H$, which has $0$ (resp. $1$) as neutral (resp. absorbing) element, and for which $\eps$ is a counit.

An \emph{interval} in $\,\EE$ is a segment in $\,\EE$ such that $(0,1):I\sqcup I\to H$ is a cofibration and $\eps:H\to I$ a weak equivalence.\end{dfn}

 The axioms of a segment are expressed by the commutativity of the following five diagrams (all isomorphisms being induced by the symmetric monoidal structure):

\begin{diagram}[small,silent](H\otimes H)\otimes H&\rTo^\sim&H\otimes(H\otimes H)\\\dTo^{\vee\otimes H}&&\dTo_{H\otimes\vee}\\H\otimes H&\overset{\vee}{\lra} H\overset{\vee}{\longleftarrow}&H\otimes H\end{diagram}

\begin{diagram}[small,silent]I\otimes H&\rTo^{0\otimes H}& H\otimes H&\lTo^{H\otimes 0}&H\otimes I\\&\rdTo_\sim&\dTo_\vee&\ldTo_\sim&\\&&H&&\end{diagram}

\begin{diagram}[small,silent]&&I\otimes H&\rTo^{1\otimes H}&H\otimes H&\lTo^{H\otimes 1}&H\otimes I&&\\&\ldTo^{I\otimes\eps}&\dTo&&\dTo_\vee&&\dTo&\rdTo^{\eps\otimes I}&\\I\otimes I&\rTo^\sim&I&\rTo^1&H&\lTo^1&I&\lTo^\sim&I\otimes I\end{diagram}

\begin{diagram}[small,silent]H\otimes H&\rTo^{\eps\otimes\eps}&I\otimes I&\quad&I&\rTo^0&H\\\dTo^\vee&&\dTo_\sim&\quad&\dTo_1&\rdTo^{id}&\dTo_\eps\\H&\rTo^\eps&I&\quad&H&\rTo^\eps&I\end{diagram}

\vspace{1ex}

\begin{exms}In many familiar contexts there is a natural choice for an interval in the sense of definition 4.1. For instance, the usual unit interval of real numbers is an interval in this sense in the monoidal model category of compactly generated spaces. The representable simplicial set $\Delta^1$ is an interval in the monoidal model category of simplicial sets. The normalised chains on $\Delta^1$ form an interval in the category of chain complexes. We will come back to these examples in more detail in Section 8.

More generally, under some general hypotheses, it can be shown that the category of segments in $\EE$ carries a Quillen closed model structure, for which the forgetful functor to $\EE$ (in fact, to $(I\sqcup I)/\EE/I$) preserves fibrations and weak equivalences, and for which every cofibrant object $H$ has the property that the map $I\sqcup I\to H$ is a cofibration in $\EE$. For example, the argument of \cite{SS} for monoids applies to segments as well, and shows that this is the case if $\EE$ satisfies the monoid axiom of loc. cit. Also, the argument of  \cite[pg. 813]{BM} shows that this is the case if $\EE$ has a monoidal fibrant replacement functor, and if in addition, the folding map $I\sqcup I\to I$ factors as a cofibration followed by a weak equivalence  $I\sqcup I\ito K\eqv I$ in $\EE$, through an object $K$ and maps in the category of cocommutative comonoids in $\EE$. Thus, there is usually an ample supply of intervals.\end{exms}

Let $T$ be a planar tree with $k$ internal edges. The planar structure induces an order on the set $E(T)$ of internal edges (for which we assume some convention to have been fixed), so that we can define$$H(T)=\bigotimes_{e\in E(T)}H.$$The group $\Aut(T)$ acts on $H(T)$ (say from the right), by canonical symmetries of $\EE$. If $D\subseteq E(T)$ is a set of internal edges, we define$$H_D(T)=\bigotimes_{e\in E(T)}H_e,$$where $H_e=I$ if $e\in D$ and $H_e=H$ if $e\not\in D$. Then the trivial cofibrations $I\overset{0}{\lra}H$ for $e\in D$ and $H\overset{id_H}{\lra} H$ for $e\not\in D$ together induce a trivial cofibration $$H_D(T)\ito H(T)$$for each subset $D\subseteq E(T)$. The pushout-product lemma, applied to the trivial cofibration $0:I\to H$, one copy for each $e\in E(T)$, then yields a trivial cofibration $$H^-(T)=\bigcup_{D\not=\emptyset}H_D(T)\ito H(T)$$Observe also that, if $T/D$ is the tree obtained by contracting the edges in $D$, there is a natural isomorphism\begin{gather}\label{collapse}H_D(T)\overset{\sim}{\lra} H(T/D)\end{gather} 

Moreover, if $c$ is a non-empty set of unary vertices in $T$, we obtain a map\begin{gather}\label{collapse2}H(T)\lra H(T/c)\end{gather}This map is defined in terms of $\vee:H\otimes H\to H$ (for vertices connecting two internal edges) and $\eps:H\to I$ (for vertices connecting an internal with an external edge), and is well-defined because $\vee$ is associative and $\eps$ is a counit for $\vee$; it depends on a convention for ordering the edges of $T$ if $\vee$ is not commutative.\vspace{1ex}

Since we are interested in the homotopical properties of the $\WW$-construction, we assume that $H$ is an \emph{interval} in the above sense, and that $P$ is a \emph{well-pointed $\Sg$-cofibrant} operad (cf. Sections \ref{cofoperad} and $3$), although for the construction itself it would be sufficient to assume that $H$ be a segment and $P$ an operad. We construct the operad $\WW(H,P)$ as a sequential colimit of trivial cofibrations of collections$$\WW_0(H,P)\ito\WW_1(H,P)\ito\WW_2(H,P)\ito\cdots$$where the filtration degree is given by the number of internal edges in the tree: $\WW_k(H,P)$ is the part of $\WW(H,P)$ which can be constructed by trees with $\leq k$ internal edges. To begin with, we set$$\WW_0(H,P)(n)=P(n)\quad\quad(n\geq 0)$$Next, suppose that $k>0$ and that $\WW_{k-1}(H,P)$ has been defined. Also suppose that this construction comes together with a canonical map\begin{gather}\label{attach}\al_S:(H(S)\otimes\underline{P}(S))\otimes_{\Aut(S)}I[\Sg_n]\lra\WW_{k-1}(H,P)(n)\end{gather}for 
each tree $S$ with $\leq k-1$ internal edges and $n$ input edges. For $k=1$, these are the unit tree $|$, and the trees $t_n$, and we take for $\al_|$ the unit $I\to P(1)$, and for $\al_{t_n}:P(n)\to P(n)$ the identity.

Consider, for a tree $T$ with $n$ input edges and $k$ internal edges, the maps of $\Aut(T)$-objects $H^-(T)\ito H(T)$ and $\underline{P}^-(T)\ito\underline{P}(T)$. The first is a trivial cofibration in $\EE$, while the second is an $\Aut(T)$-cofibration, see Section 3(\ref{minus2}). Consequently, by Lemma \ref{eq2}, the map of the pushout-product axiom\begin{gather}\label{glue}(H\otimes P)^-(T)\ito H(T)\otimes\underline{P}(T)\end{gather} is a trivial $\Aut(T)$-cofibration, where for brevity we have written$$(H\otimes P)^-(T)=(H^-(T)\otimes\underline{P}(T))\cup_{H^-(T)\otimes \underline{P}^-(T)}(H(T)\otimes\underline{P}^-(T))$$Next, we show that the maps $\al_S$ in (\ref{attach}) together define a $\Sg_n$-equivariant map\begin{gather}\label{attach2}\al_T^-:(H\otimes P)^-(T)\otimes_{\Aut(T)}I[\Sg_n]\lra\WW_{k-1}(H,P)(n)\end{gather}Indeed, for each non-empty subset $D$ of internal edges, the map $H_D(T)\overset{\sim}{\lra}H(T/D)$ in (\ref{collapse}) and the operad composition map $\underline{P}(T)\lra\underline{P}(T/D)$, together with the map $\al_{T/D}$, constructed inductively in (\ref{attach}), define a map$$(H_D(T)\otimes\underline{P}(T))\otimes_{\Aut(T)}I[\Sg_n]\lra\WW_{k-1}(H,P)(n)$$These maps together give a well-defined $\Sg_n$-equivariant map\begin{gather}\label{attach3}(H^-(T)\otimes\underline{P}(T))\otimes_{\Aut(T)}I[\Sg_n]\lra\WW_{k-1}(H,P)(n)\end{gather}Similarly, for each nonempty set $c$ of unary vertices in $T$, the maps $\underline{P}_c(T)\to\underline{P}(T/c)$ and $H(T)\to H(T/c)$ in (\ref{collapse2}), together with $\al_{T/c}$ define a map$$(H(T)\otimes\underline{P}_c(T))\otimes_{\Aut(T)}I[\Sg_n]\lra\WW_{k-1}(H,P)(n)$$These maps together yield a well-defined $\Sg_n$-equivariant map\begin{gather}\label{attach4}(H(T)\otimes\underline{P}^-(T))\otimes_{\Aut(T)}I[\Sg_n]\lra\WW_{k-1}(H,P)(n)\end{gather}The maps (\ref{attach3}) and (\ref{attach4}) together give the map (\ref{attach2}). Now, take the coproduct over isomorphism classes of such trees $T$ (with $n$ input edges and $k$ internal edges), and the coproduct of the trivial $\Aut(T)$-cofibrations in (\ref{glue}) and the coproduct of the maps of form (\ref{attach2}), and construct the pushout\begin{gather}\label{W}\begin{diagram}[small,silent]\coprod_{[T],T\in\TT(n,k)}(H\otimes P)^-(T)\otimes_{\Aut(T)}I[\Sg_n]&\rTo^{\coprod\al^-_T}&\WW_{k-1}(H,P)(n)\\\dTo&&\dTo\\\coprod_{[T],T\in\TT(n,k)}(H(T)\otimes\underline{P}(T))\otimes_{\Aut(T)}I[\Sg_n]&\rTo^{\coprod\al_T}&\WW_k(H,P)(n)\end{diagram}\end{gather}This defines $\WW_k(H,P)(n)$, as well as the maps $\al_T$ as the restrictions of the bottom map in (\ref{W}). Note that by construction, $\WW_{k-1}(H,P)(n)\ito\WW_k(H,P)(n)$ is a trivial $\Sg_n$-cofibration. This completes the definition of the sequence$$\WW_0(H,P)(n)\ito\WW_1(H,P)(n)\ito\WW_2(H,P)(n)\ito\cdots$$and hence of the colimit $\WW(H,P)(n)=\varinjlim\WW_k(H,P)(n)$. By construction, each inclusion $\WW_k(H,P)(n)\ito\WW(H,P)(n)$ is again a trivial $\Sg_n$-cofibration. In particular, for $k=0$, we find that$$P(n)\ito\WW(H,P)(n)$$ is a trivial $\Sg_n$-cofibration.

Finally, we observe that there is an operad structure on $\WW(H,P)$. One way to describe it is in terms of the maps$$\al_T:(H(T)\otimes\underline{P}(T))\otimes_{\Aut(T)}I[\Sg_n]\lra\WW(H,P)(n)$$Indeed, given a tree $T$ with $n$ input edges and $n$ trees $T_1,\dots,T_n$, with $k_i$ input edges respectively, one obtains a new tree $T'=T(T_1,\dots,T_n)$ with $k=k_1+\cdots+k_n$ input edges. The $n$ edges of $T$ become internal edges of $T'$. Assigning length $1$ to these internal edges, we obtain a map$$H(T)\otimes H(T_1)\otimes\cdots\otimes H(T_n)\lra H(T')$$Also the free operad structure on $\underline{P}$ gives a map$$\underline{P}(T)\otimes\underline{P}(T_1)\otimes\cdots\otimes\underline{P}(T_n)\lra\underline{P}(T')$$In this way we obtain a map$$(H(T)\otimes\underline{P}(T))_{\lambda[T]}\otimes\bigotimes_{i=1}^{n}(H(T_i)\otimes\underline{P}(T_i))_{\lambda[T_i]}\overset{\gamma}{\lra}(H(T')\otimes\underline{P}(T'))_{\lambda[T']}$$satisfying suitable equivariance conditions; here and below, the subscript $\lambda[T]$ stands for the tensor product $(-)\otimes_{\Aut(T)}I[\Sg_n]$. The operad structure on $\WW(H,P)$ is uniquely determined by the requirement that for each grafting operation $T'=T(T_1,\dots,T_k)$ the following diagram commutes:\begin{gather}\label{operad}\begin{diagram}[small,silent](H(T)\otimes\underline{P}(T))_{\lambda[T]}\otimes\bigotimes_{i=1}^{n}(H(T_i)\otimes\underline{P}(T_i))_{\lambda[T_i]}&\rTo&\WW(H,P)(n)\otimes\bigotimes_{i=1}^n\WW(H,P)(k_i)\\\dTo^\gamma&&\dTo\\(H(T')\otimes\underline{P}(T'))_{\lambda[T']}&\rTo^{\al_{T'}}&\WW(H,P)(k)\end{diagram}\end{gather}This concludes the definition of the operad $\WW(H,P)$.

\begin{rmk}\label{filter}A \emph{filtered} operad $Q$ is an operad which is union of subcollections$$Q^{(1)}\subseteq Q^{(2)}\subseteq\cdots$$with the property that the operad multiplication restricts to maps$$Q^{(k)}(n)\otimes Q^{(l_1)}(m_1)\otimes\cdots\otimes Q^{(l_n)}(m_n)\lra Q^{(k+l_1+\cdots+l_n)}(m_1+\cdots+m_n).$$$\WW(H,P)$ becomes a filtered operad if we shift the grading, as in$$\WW(H,P)^{(k)}=\WW_{k-1}(H,P).$$This shifted filtration degree corresponds to the number of vertices of the indexing trees, whereas the original filtration degree corresponds to the number of internal edges of the indexing trees.\end{rmk}

\begin{rmk}\label{analog}We have treated the $\WW$-construction for symmetric operads without any constraint. There are other (simpler) types of operads which present some interest in their own, and for which a modified $\WW$-construction is available:\vspace{1ex}

(a) \emph{Non-symmetric operads}, i.e. operads without symmetric group actions. Here, the free operad and $\WW$-constructions simplify considerably, since there are no tree-automorphisms to be taken care of. More precisely, the forgetful functor from symmetric operads to non-symmetric operads has a left adjoint which identifies the category of non-symmetric operads with a full coreflective subcategory of the category of symmetric operads. There is an essentially unique $\WW$-construction for non-symmetric operads which is compatible with the $\WW$-construction for symmetric operads under this embedding.\vspace{1ex}


(b) \emph{Reduced operads}, i.e. operads $P$ such that $P(0)=I$. These are operads with a unique nullary operation. Often, the category of reduced operads has a model structure even if the category of all symmetric operads does not, see \cite[Theorem 3.1]{BM}. This is in particular the case for chain operads. Algebras over reduced operads behave like \emph{pointed} algebras over operads \emph{without} nullary operation. This is the reason for which reduced operads in \emph{pointed} symmetric monoidal categories are often replaced by augmented operads without nullary operation. This simplifies the corresponding free operad and $\WW$-constructions, but forces the degeneracy maps (see below) to be zero. For sake of completeness we treat the free reduced operad and reduced $\WW$-constructions in full generality.

Any reduced operad comes equipped with \emph{degeneracy maps} defined by $P(n)=P(n)\otimes I^{\otimes n}\to P(n)\otimes P(\eps_1)\otimes\cdots P(\eps_n)\to P(\eps_1+\cdots+\eps_n)$ for $\eps_k\in\{0,1\},1\leq k\leq n$, cf. May \cite{M}. These degeneracy maps together with the symmetric group actions assemble into a contravariant functor $\Lambda_0^{op}\to\EE:\nn\mapsto P(n)$, where $\Lambda_0$ is the category of finite sets $\nn=\{1,\dots,n\},n\geq 0,$ and injections; such a contravariant functor has been called a \emph{pre-operad} by the first named author \cite{B}. 

The category of reduced operads in $\EE$ is then monadic over the category of pre-operads in $\EE$. In this setting, the free reduced operad and reduced $\WW$-constructions are literally the same as the constructions presented in Sections $3$ and $4$, except that all trees with input edges and some vertices of valence $0$ have to be discarded; instead, the given pre-operad structure extends in a natural way to both, the free reduced operad and reduced $\WW$-constructions. Boardman and Vogt \cite[Chapter V.3]{BV} also consider reduced versions of their topological $\WW$-construction. Their $\WW^{''}$-contruction coincides for $\EE=\Top$ with the reduced $\WW$-construction just described.\vspace{1ex}   

(c) \emph{Pseudo-operads}, i.e. ``operads'' without nullary operations and without unit, cf. \cite[Sections 1.3 and 1.7]{MSS}.  The free operad and $\WW$-constructions for pseudo-operads are similar to the corresponding constructions for general operads: one has to discard the unit-tree $|$, all trees with vertices of valence $0$ and all unit identifications, cf. \cite[Section 1.9]{MSS}. Therefore, the binary operation of the segment is no longer needed to carry out the $\WW$-construction for pseudo-operads.\end{rmk}

\section{Cofibrancy of the $\WW$-construction}

Suppose $\EE$ is a monoidal model category with interval $H$. We have just constructed an operad $\WW(H,P)$ for each well-pointed $\Sg$-cofibrant operad $P$. In this section we show that $\WW(H,P)$ provides a cofibrant resolution for $P$. The following stronger statement is given in terms of the classes of weak equivalences, fibrations and cofibrations which we defined in Section \ref{cofoperad} for the category of operads in $\EE$ (even if these do not always form a model structure).

\begin{thm}\label{main}Let $\EE$ be a cofibrantly generated monoidal model category with cofibrant unit $I$ and interval $H$. For any well-pointed $\Sg$-cofibrant operad $P$, the counit of the adjunction between pointed collections and operads admits a factorisation$$\FF_*(P)\overset{\delta}{\ito}\WW(H,P)\overset{\gamma}{\lra}P$$ into a cofibration $\delta$ followed by a weak equivalence $\gamma$. In particular, $\WW(H,P)$ is a cofibrant resolution for $P$.\end{thm}

\begin{proof}The second assertion follows from the first by Lemma \ref{cofibrant}. Let us begin by describing the map $\gamma$ explicitly. The counit of the interval $H$ and the operad structure of $P$ define a map $H(T)\otimes\underline{P}(T)\to P(n)$ for each tree $T$ with $n$ input edges. These fit together into a map of operads $\gamma:\WW(H,P)\to P$ such that the composition $\WW_0(H,P)\ito\WW(H,P)\overset{\gamma}{\lra}P$ is the identity. Since $\WW_0(H,P)\ito\WW(H,P)$ is a trivial cofibration of collections, $\gamma$ is a weak equivalence.

Next, the inclusions $\underline{P}(T)\ito H(T)\otimes\underline{P}(T)$ induced by $1:I\to H$ define a map of operads $\delta:\FF_*(P)\to\WW(H,P)$. (Intuitively, $\delta$ views the elements of $\FF_*(P)$ as belonging to $\WW(H,P)$ by giving all internal edges length $1$). Clearly, the composition $\gamma\delta$ may be identified with the counit of the adjunction between pointed collections and operads, so it remains to be shown that $\delta$ is a cofibration of operads. In other words, we have to show that in any commutative square of operad maps\begin{diagram}[small,silent]\FF_*(P)&\rTo&Y\\\dTo^\delta&&\dTo_\chi\\\WW(H,P)&\rTo^\psi&X\end{diagram}where $\chi$ is a trivial fibration, a diagonal filler exists. By adjunction, an operad map $\FF_*(P)\to Z$ corresponds to a map of collections $\WW_0(H,P)=P\to Z$. Thus, the square of operad maps corresponds to a square of maps of collections:\begin{diagram}[small,silent]\WW_0(H,P)&\rTo^{\phi_0}&Y\\\dTo&&\dTo_{\chi}\\\WW(H,P)&\rTo^{\psi}&X\end{diagram}An operad map $\WW(H,P)\to Y$ is a filler for the former square if and only if the underlying map of collections is a filler for the latter. The existence of such a filler follows  from Lemma \ref{khom}, by induction on the filtration degree of $\WW(H,P)$.\end{proof}For the proof of Lemma \ref{khom}, we introduce the concept of a \emph{$k$-homomorphism}. Recall that $\WW(H,P)$ is filtered by subcollections $\WW_k(H,P)$. Informally, a map of collections $\phi_k:\WW_k(H,P)\to Y$ into any operad $Y$ is a $k$-homomorphism if internal edges of length $1$ in $\WW_k(H,P)$ are sent to operad compositions in $Y$. For a formal definition, the following notation will be useful:

For a pair of trees $T_1$ and $T_2$ with $n_1$ and $n_2$ input edges respectively, and $k_1$ and $k_2$ internal edges, build a new tree $T=T_1\circ_eT_2$ by grafting $T_2$ onto the input edge $e$ of $T_1$. The tree $T$ will have $k=k_1+k_2+1$ internal edges and $n=n_1+n_2-1$ input edges. There is a map $H(T_1)\otimes H(T_2)\to H(T)$ giving the new internal edge $e$ of $T$ length $1$. Also, there is an obvious map $\underline{P}(T_1)\otimes\underline{P}(T_2)\to\underline{P}(T)$. Together, these define a map$$\bigotimes_{i=1,2}(H(T_i)\otimes\underline{P}(T_i))\otimes_{\Aut(T_i)}I[\Sg_{n_i}]\overset{\gamma_e}{\lra}(H(T)\otimes\underline{P}(T))\otimes_{\Aut(T)}I[\Sg_n]$$Also, for $j=1,\dots,n_1$, there are operations $Y(n_1)\otimes Y(n_2)\overset{\circ_j}{\lra}Y(n)$ induced by the operad structure of $Y$, and these give a map$$\coprod_{j=1,\dots,n_1}Y(n_1)\otimes Y(n_2)\overset{\circ}{\lra}Y(n)$$

Now we can define the notion of a \emph{$k$-homomorphism} into $Y$ by induction on $k$: 

\begin{dfn}Let $Y$ be an operad. A \emph{$0$-homomorphism} is a map of pointed collections $\phi_0:\WW_0(H,P)\to Y$. For $k>0$, a \emph{$k$-homomorphism} is a map of collections $\phi_k:\WW_k(H,P)\to Y$ such that\vspace{1ex}

(a) for each $l<k$, the restriction of $\phi_k$ along $\WW_l(H,P)\ito\WW_k(H,P)$ is an $l$-homomorphism;

(b) for each tree $T=T_1\circ_eT_2$ as above, the following diagram commutes:\begin{diagram}\bigotimes_{i=1,2}(H(T_i)\otimes\underline{P}(T_i))\otimes_{\Aut(T_i)}I[\Sg_{n_i}]&\rTo^{\gamma_e}&(H(T)\otimes\underline{P}(T))\otimes_{\Aut(T)}I[\Sg_n]\\\dTo_{\al_{T_1}\otimes\al_{T_2}}&&\dTo_{\al_T}\\\coprod_{j=1}^{n_1}\WW_{k_1}(H,P)\otimes\WW_{k_2}(H,P)&&\WW_k(H,P)\\\dTo_{\phi_{k_1}\otimes\phi_{k_2}}&&\dTo_{\phi_k}\\\coprod_{j=1}^{n_1}Y(n_1)\otimes Y(n_2)&\rTo^{\circ}&Y\end{diagram}\end{dfn}

Here, the maps $\al_{T_i}$ and $\al_T$ are as described in Section $4$. Notice that $(H(T_1)\otimes\underline{P}(T_1))\otimes_{\Aut(T_1)}I[\Sg_{n_1}]$ decomposes into a sum of $n_1$ factors, each factor being determined by the number $j$ assigned to the input edge $e$ under the identification of $I[\Sg_{n_1}]$ with $I[\lambda(T_1)]$. The upper left map in the diagram sends this summand to the summand indexed by the same $j$. Notice that a sequence of $k$-homomorphisms $\phi_k:\WW_k(H,P)\to Y$ such that for each $l<k$, the map $\phi_l$ equals the restriction$\WW_l(H,P)\ito\WW_k(H,P)\overset{\phi_k}{\lra}Y$, determines a unique operad map $\WW(H,P)\to Y$ in the colimit. In fact, a map of collections $\WW(H,P)\to Y$ is an operad map if and only if each of its restrictions $\WW_k(H,P)\to Y$ is a $k$-homomorphism.

\begin{rmk}The suboperads of $\WW(H,K)$ generated by the $\WW_k(H,P)$ also admit a nice description in terms of trees, and could be used to give an alternative proof of Theorem \ref{main}, cf. the proof of \cite[Theorem 4.1]{V}.\end{rmk} 

\begin{lma}\label{khom}Let $k> 0$ and let\begin{diagram}[small,silent]\WW_{k-1}(H,P)&\rTo^{\phi_{k-1}}&Y\\\dTo&&\dTo_\chi\\\WW(H,P)&\rTo^\psi&X\end{diagram} be a commutative square, where $\psi$ is a map of operads, $\chi$ is a trivial fibration of operads, and $\phi_{k-1}$ is a $(k-1)$-homomorphism. Then there exists a $k$-homomorphism $\phi_k$ extending $\phi_{k-1}$ such that  the following diagram commutes:\begin{diagram}[small,silent]\WW_{k-1}(H,P)&\rTo&\WW_k(H,P)&\rTo^{\phi_k}&Y\\&\rdTo&\dTo&&\dTo_\chi\\&&\WW(H,P)&\rTo^\psi&X\end{diagram}\end{lma}

\begin{proof}Consider a tree $T$ with $n$ input edges and $k$ internal vertices. Recall the cofibrations $H^-(T)\ito H(T)$ and $H_D(T)\ito H(T)$ for each set $D$ of internal edges in $T$. For each internal edge $e\in E(T)$, we already defined$$H_e=\begin{cases}I&\text{if }e\in D\\H&\text{if }e\not\in D\end{cases}$$and now set$$H^+_e=\begin{cases}I\sqcup I&\text{if }e\in D\\H&\text{if }e\not\in D\end{cases}$$Then $H_D(T)=\bigotimes_eH_e, H(T)=\bigotimes_eH$ as before, and we define $$H_D^+(T)=\bigotimes_eH^+_e$$Thus $I\overset{0}{\ito}I\sqcup I\overset{(0,1)}{\ito}H$ define cofibrations$$H_D(T)\ito H^+_D(T)\ito H(T)$$Also, whereas $H^-(T)=\bigcup_{D\not=\emptyset}H_D(T)$ as before, we define$$H^+(T)=\bigcup_{D\not=\emptyset}H^+_D(T)$$Then the pushout-product lemma gives cofibrations$$H^-(T)\ito H^+(T)\ito H(T)$$ and these are $\Aut(T)$-equivariant maps. Define$$(H\otimes P)^+(T)=(H^+(T)\otimes\underline{P}(T))\cup_{(H^+(T)\otimes\underline{P}^-(T))}(H(T)\otimes\underline{P}^-(T))$$so that we have $\Aut(T)$-cofibrations$$(H\otimes P)^-(T)\ito(H\otimes P)^+(T)\ito H(T)\otimes\underline{P}(T)$$Now consider the pushouts of Section $4$ which defined $\WW_k(H,P)(n)$ out of the preceding filtration term $\WW_{k-1}(H,P)(n)$. There is a factorisation of the $\Sg_n$-cofibration $\WW_{k-1}(H,P)(n)\ito\WW_k(H,P)(n)$ by the pushout

\begin{diagram}[small,silent]\coprod_{[T]}(H\otimes P)^-(T)\otimes_{\Aut(T)}I[\Sg_n]&\rTo^{\coprod\al_T^-}&\WW_{k-1}(H,P)(n)\\\dTo&&\dTo\\\coprod_{[T]}(H\otimes P)^+(T)\otimes_{\Aut(T)}I[\Sg_n]&\rTo^{\coprod\al_T^+}&\WW^+_{k-1}(H,P)(n)\\\dTo&&\dTo\\\coprod_{[T]}(H(T)\otimes\underline{P}(T))\otimes_{\Aut(T)}I[\Sg_n]&\rTo^{\coprod\al_T}&\WW_k(H,P)(n)\end{diagram}This defines a refinement of the filtration of $\WW(H,P)$ by subcollections $\WW_k^+(H,P)$ and further (non-trivial) cofibrations\begin{diagram}[small,silent]\WW_0(H,P)&\rTo&\WW_1(H,P)&\rTo&\WW_2(H,P)&\rTo&\WW_3(H,P)\,\cdots\\\dTo&\ruTo&\dTo&\ruTo&\dTo&\ruTo\\\WW_0^+(H,P)&&\WW_1^+(H,P)&&\WW^+_2(H,P)&\cdots\end{diagram}Clearly, any $(k-1)$-homomorphism $\phi_{k-1}:\WW_{k-1}(H,P)\to Y$ extends \emph{uniquely} to a map $\phi_{k-1}^+:\WW_{k-1}^+(H,P)\to Y$ which satisfies the conditions for a $k$-homomorphism along the maps $\al^+_T$. Since $\chi:Y\to X$ is a map of operads, this extension makes the following diagram commute:\begin{diagram}[small,silent]\WW^+_{k-1}(H,P)&\rTo^{\phi^+_{k-1}}&Y\\\dTo&&\dTo_\chi\\\WW_k(H,P)&\rTo&X\end{diagram}Now take any diagonal filler $\phi_k$ in this diagram, which exists by the model structure on collections. This is the required $k$-homomorphism.\end{proof}

\section{Functoriality of the $\WW$-construction}

The $\WW$-construction is a functor in two variables: the segment and the operad. In this section, we study the homotopical behaviour of the $\WW$-construction with respect to maps of segments $H\to K$, with respect to maps of operads $P\to Q$, and with respect to symmetric monoidal functors $\EE\to\EE'$. We assume throughout that our ambient category $\EE$ (resp. $\EE'$) is a cofibrantly generated monoidal model category with cofibrant unit $I$ (resp. $I'$).\vspace{1ex}

It will be convenient to call a segment $H$ \emph{cofibrant} whenever $I\sqcup I\to H$ is a cofibration. The underlying object of a cofibrant segment is cofibrant.

Observe that the unit $I$ is a segment $I\sqcup I\overset{\nabla}{\lra}I\overset{id}{\lra}I$ with binary operation given by the canonical isomorphism $I\otimes I\to I$.

The coproduct $I\sqcup I$ is a cofibrant segment $I\sqcup I\overset{id}{\lra}I\sqcup I\overset{\nabla}{\lra}I$ with binary operation induced by the maximum operation on the two-element set $\{0,1\}$ via the isomorphism $I[\{0,1\}]\cong I\sqcup I$. The codiagonal $I\sqcup I\to I$ is then a map of segments.

\begin{lma}\label{counit}For each operad $P$, the codiagonal $I\sqcup I\to I$ induces a map of operads $\WW(I\sqcup I,P)\to\WW(I,P)$ which is canonically isomorphic to the counit $\FF_*(P)\to P$ of the adjunction between pointed collections and operads.\end{lma}

\begin{proof}With the notation of the proof of Lemma \ref{khom}, we have for each $k>0$, equality $\WW^+_{k-1}(I\sqcup I,P)=\WW_k(I\sqcup I,P)$. This implies that for each operad $Y$ and each $0$-homomorphism $\phi_0:\WW_0(I\sqcup I,P)\to Y$, there exist uniquely determined, mutually compatible $k$-homomorphisms $\WW_k(I\sqcup I,P)\to Y$ extending $\phi_0$. Therefore, there exists a unique map of operads $\phi:\WW(I\sqcup I,P)\to Y$ extending $\phi_0$. This is precisely the universal property characterising $\FF_*(P)$, cf. Section 3. 

The counit $\FF_*(P)\to P$ corresponds under adjunction to the identity of the pointed collection underlying $P$. If we put $Y=P$ and $\phi_0=id_P$, the above extension yields the augmentation $W(I\sqcup I,P)\to P$. The latter is easily identified with the map of operads $W(I\sqcup I,P)\to\WW(I,P)$ induced by the codiagonal $I\sqcup I\to I$.\end{proof}

In the study of the Godement resolution in Section \ref{Godement0}, we need an extension of the preceding lemma, based on the following definition:
\begin{dfn}For any segment $H$, let $H^\diamond=H\sqcup I$ be the segment $$I\sqcup I\overset{0_H\sqcup id_I}{\lra}H^\diamond\overset{(\eps_H,id_I)}{\lra}I$$ with binary operation $H^\diamond\otimes H^\diamond\to H^\diamond$ induced by the binary operation of $H$ on $H\otimes H$, by the counit of $H$ on $H\otimes I$ and $I\otimes H$, and by the canonical isomorphism $I\otimes I\cong I$ on the last component of $H^\diamond\otimes H^\diamond$.\end{dfn}Observe that $H^\diamond$ is obtained from $H$ by adding an external absorbing element; it is a cofibrant segment whenever $H$ is. In particular, $I^\diamond$ is precisely the segment $I\sqcup I$. Moreover, for each segment $H$, there is map of segments $(id_H,1_H):H^\diamond\to H$.

\begin{lma}\label{diamond}For any segment $H$ and operad $P$, the map  $\WW(H^\diamond,P)\to\WW(H,P)$ induced by $H^\diamond\to H$ is canonically isomorphic to the counit $\FF_*(\WW(H,P))\to\WW(H,P)$ at $\WW(H,P)$ of the adjunction between pointed collections and operads.\end{lma}

\begin{proof}The canonical cofibration $H\to H^\diamond$ preserves $0$ and the binary operation, and induces thus a map of pointed collections $\WW(H,P)\to\WW(H^\diamond,P)$. It suffices to show that any map of pointed collections $\WW(H,P)\to Y$ into an operad $Y$ uniquely extends along $\WW(H,P)\to\WW(H^\diamond,P)$ to a map of operads $\WW(H^\diamond,P)\to Y$. This extension is constructed like in the proof of Lemma \ref{counit} by means of a compatible family of $k$-homomorphisms $\WW_k(H^\diamond,P)\to Y$ extending the given $\WW_k(H,P)\to Y$. In order to construct this family we proceed in several steps. It follows from the definitions that for each tree $T$, the left square below induces the two squares on the right, where free use is made of the notations of Section 5:\begin{diagram}[small,silent]I&\rTo^0&H&\quad\quad&H^-(T)&\rTo&H(T)&\quad\quad&(H\otimes P)^-(T)&\rTo&H(T)\otimes\underline{P}(T)\\\dTo^0&&\dTo&&\dTo&&\dTo&&\dTo&&\dTo\\I\sqcup I&\rTo&H^\diamond&&H^{\diamond+}(T)&\rTo&H^\diamond(T)&&(H^\diamond\otimes P)^+(T)&\rTo&H^\diamond(T)\otimes\underline{P}(T)\end{diagram}The left hand square is a pushout by definition of $H^\diamond$. Since the tensor product commutes with  coproducts on both sides, the middle square is also a pushout. The right square is obtained from the middle square by tensoring with $\underline{P}^-(T)\to\underline{P}(T)$ and taking pushout-product maps; it is therefore a pushout by Lemma \ref{cube} below. Diagram (\ref{W}) of Section 4 then yields the following pushout diagram of collections\begin{diagram}[small,silent]\WW_{k-1}(H,P)&\rTo&\WW_k(H,P)\\\dTo&&\dTo\\\WW_{k-1}^+(H^\diamond,P)&\rTo&\WW_k(H^\diamond,P)\end{diagram}from which the existence and unicity of the required family of $k$-homomorphisms immediately follows. Extending in this way the identity map of $\WW(H,P)$ yields the map $\WW(H^\diamond,P)\to\WW(H,P)$ induced by $H^\diamond\to H$.\end{proof}

\begin{lma}\label{invariance}For any map of intervals $H\to K$ and any weak equivalence of well-pointed $\Sg$-cofibrant operads $P\to Q$, the induced map $\WW(H,P)\to\WW(K,Q)$ is a weak equivalence of cofibrant operads.\end{lma}

\begin{proof}Apply Theorem \ref{main} and the 2-out-of-3 property of weak equivalences.\end{proof}

\begin{prp}\label{cofibrant1}For any (trivial) cofibration of cofibrant segments $H\to K$ and any well-pointed $\Sg$-cofibrant operad $P$, the induced map $\WW(H,P)\to\WW(K,P)$ is a (trivial) cofibration of operads.\end{prp}

\begin{proof}We have to show that $\WW(H,P)\to\WW(K,P)$ has the left lifting property with respect to any trivial fibration (resp. fibration) of operads $\chi:Y\to X$. The required diagonal filler is constructed by means of a compatible family of $k$-homomorphisms. To be precise, we consider the following commutative diagram:\begin{diagram}[small,silent]\WW_{k-1}(H,P)&\rTo&\WW_{k-1}^+(H,P)&\rTo^{i_k^H}&\WW_k(H,P)&\rTo&Y\\\dTo&&\dTo&&\dTo&&\dTo_\chi\\\WW_{k-1}(K,P)&\rTo&\WW_{k-1}^+(K,P)&\rTo^{i_k^K}&\WW_k(K,P)&\rTo&X\end{diagram}Since $\WW_0(H,P)=P=\WW_0(K,P)$ we may assume by induction that there is a $(k-1)$-homomorphism $\phi_{k-1}:\WW_{k-1}(K,P)\to Y$ making the two triangles forming the outer rectangle commute. It follows from the definitions that there is a \emph{unique} extension of $\phi_{k-1}$ to a $k$-homomorphism $\phi_{k-1}^+:\WW_{k-1}^+(K,P)\to Y$. In order to extend $\phi_{k-1}^+$ to the required $k$-homomorphism $\phi_k:\WW_k(K,P)\to Y$, it suffices to show that 
\begin{gather}\label{coprod}\WW_{k-1}^+(K,P)\cup_{\WW^+_{k-1}(H,P)}\WW_k(H,P)\to\WW_k(K,P)\end{gather}is a (trivial) cofibration of collections. As we have seen in the proof of Lemma \ref{khom}, the maps $i_k^H$ and $i_k^K$ are pushouts, which depend functorially on $K\to H$. Lemma \ref{cube} then implies that (\ref{coprod}) is part of the following pushout diagram of collections:\begin{gather}\label{chase}\begin{diagram}[small,silent]\coprod_{[T]}L(T)\otimes_{\Aut(T)}I[\Sg_{|T|}]&\rTo&\WW_{k-1}^+(K,P)\cup_{\WW^+_{k-1}(H,P)}\WW_k(H,P)\\\dTo&&\dTo\\\coprod_{[T]}(K(T)\otimes\underline{P}(T))\otimes_{\Aut(T)}I[\Sg_{|T|}]&\rTo&\WW_k(K,P)\end{diagram}\end{gather}The left vertical map of (\ref{chase}) is obtained from the commutative square\begin{diagram}[small,silent](H\otimes P)^+(T)&\rTo&(K\otimes P)^+(T)\\\dTo&&\dTo\\H(T)\otimes \underline{P}(T)&\rTo&K(T)\otimes\underline{P}(T)\end{diagram}by taking the pushout comparison map $L(T)\to  K(T)\otimes\underline{P}(T)$, applying the left Quillen functor $(-)\otimes_{\Aut(T)}I[\Sg_{|T|}]$ and summing up over isomorphism classes of trees with $k$ internal edges. It will thus be sufficient to show that for each tree $T$, the map $L(T)\to  K(T)\otimes\underline{P}(T)$ is a (trivial) $\Aut(T)$-cofibration. Observe first that this map is also the pushout-product map of the following commutative square:\begin{diagram}[small,silent](H(T)\cup_{H^+(T)}K^+(T))\otimes\underline{P}^-(T)&\rTo&(H(T)\cup_{H^+(T)}K^+(T))\otimes\underline{P}(T)\\\dTo&&\dTo\\K(T)\otimes\underline{P}^-(T) &\rTo&K(T)\otimes\underline{P}(T)\end{diagram}
For well-pointed $\Sg$-cofibrant operads $P$, the map $\underline{P}^-(T)\to\underline{P}(T)$ is an $\Aut(T)$-cofibration. Therefore, by Lemma \ref{eq2}, it remains to be shown that $H(T)\cup_{H^+(T)}K^+(T)\to K(T)$ is a (trivial) cofibration. This follows from Lemma \ref{eq4}, with $i$ running through the set of internal edges of $T$, $\Gamma$ being the trivial group, and $A_i\ito B_i\ito C_i$ being $I\sqcup I\ito H\ito K$.\end{proof} 

\begin{prp}\label{cofibrant2}For any cofibrant segment $H$ and any (trivial) $\Sg$-cofibration of well-pointed $\Sg$-cofibrant operads $P\to Q$, the induced map $\WW(H,P)\to\WW(H,Q)$ is a (trivial) cofibration of operads.\end{prp}

\begin{proof} By the same argument as in the proof of Proposition \ref{cofibrant1}, it is sufficient to show that for each tree $T$, the pushout comparison map of\begin{diagram}[small,silent](H\otimes P)^+(T)&\rTo&(H\otimes Q)^+(T)\\\dTo&&\dTo\\H(T)\otimes \underline{P}(T)&\rTo&H(T)\otimes\underline{Q}(T)\end{diagram} is a (trivial) $\Aut(T)$-cofibration. Again, this map is also the pushout-product map of the following commutative square:

\begin{diagram}[small,silent]H^+(T)\otimes(\underline{P}(T)\cup_{\underline{P}^-(T)}\underline{Q}^-(T))&\rTo&H^+(T)\otimes\underline{Q}(T)\\\dTo&&\dTo\\H(T)\otimes(\underline{P}(T)\cup_{\underline{P}^-(T)}\underline{Q}^-(T)) &\rTo&H(T)\otimes\underline{Q}(T)\end{diagram}For cofibrant segments $H$, the map $H^+(T)\to H(T)$ is a cofibration. Therefore, by Lemma \ref{eq2}, it remains to be shown that for each tree $T$, the map $$\phi_T:\underline{P}(T)\cup_{\underline{P}^-(T)}\underline{Q}^-(T)\to\underline{Q}(T)$$ is a (trivial) $\Aut(T)$-cofibration. We use induction on trees: For the unit-tree $|$, we get $\phi_{|}=id_I$; for the tree $t_n$ with one vertex and $n$ input edges, $\phi_{t_n}$ is the given (trivial) $\Sg_n$-cofibration $P(n)\to Q(n)$. Now, let $T$ be a tree obtained by grafting $T=t_n(T_1,\dots,T_n)$. If $n=1$, then $\phi_T$ is the map$$(P(1)\otimes\underline{ P}(T_1))\cup (I\otimes\underline{Q}(T_1))\cup (Q(1)\otimes\underline{Q}^-(T_1))\to Q(1)\otimes\underline{Q}(T_1)$$which is a (trivial) $\Aut(T)$-cofibration by Lemma \ref{eq4}, putting $\Gamma=\Aut(T_1)=\Aut(T),$ and letting $A_1\ito B_1\ito C_1$ be the cofibrations $I\ito P(1)\ito Q(1)$ (with trivial $\Gamma$-action), and $A_2\ito B_2\ito C_2$ be the $\Gamma$-cofibrations $\underline{Q}^-(T_1)\ito\underline{P}(T_1)\cup_{\underline{P}^-(T_1)}\underline{Q}^-(T_1)\ito Q(T_1)$. If $n>1$, then $\phi_T$ is the map$$(P(n)\otimes\underline{P}(T^{up}))\cup (P(n)\otimes\underline{Q}^-(T^{up}))\cup (Q(n)\otimes \underline{Q}^-(T^{up}))\to Q(n)\otimes \underline{Q}(T^{up})$$where we use for any collection $\KK$ the notations $\underline{\KK}(T)=\KK(n)\otimes\underline{\KK}(T^{up})$ and $\underline{\KK}^-(T)=\KK(n)\otimes\underline{\KK}^-(T^{up})$. In particular, the semi-direct product $\Aut(T)=\Gamma_T\rtimes\Sg_T$ of Section 3 acts in a canonical way on the tensor products above. An inductive application of Lemma \ref{eq1} and Lemma \ref{eq5} with respect to the split exact sequence $1\to\Gamma_T\to\Gamma_T\rtimes\Sg_T\to\Sg_T\to 1$ implies then that $\phi_T$ is a (trivial) $\Aut(T)$-cofibration as required, letting the $\Sg_T$-cofibrations be $0\ito P(n)\ito Q(n)$, and the $\Gamma_T$-cofibrations be $\underline{Q}^-(T^{up})\ito\underline{P}(T^{up})\cup_{\underline{P}^-(T^{up})}\underline{Q}^-(T^{up})\ito \underline{Q}(T^{up})$. \end{proof}

\begin{rmk}The two preceding propositions generalise in two different ways the fact (cf. Theorem \ref{main}) that $\WW(H,P)$ is a cofibrant operad for cofibrant segments $H$ and well-pointed $\Sg$-cofibrant operads $P$. Indeed, the cofibrancy of $\WW(H,P)$ follows either from Proposition \ref{cofibrant1} and Lemma \ref{counit}, keeping $P$ fixed and using that $I\sqcup I\to H$ is a cofibration of cofibrant segments; or, it follows from Proposition \ref{cofibrant2}, keeping $H$ fixed and using that $\underline{I}\to P$ is a $\Sg$-cofibration of well-pointed $\Sg$-cofibrant operads, and that $\WW(H,\underline{I})\cong\underline{I}$, where $\underline{I}$ denotes the initial operad.\end{rmk}

\begin{prp}\label{Quillen}Let $\Phi:\EE\to\EE'$ be a unit- and colimit-preserving symmetric monoidal functor. For each operad $P$ and each segment $H$ in $\,\EE$, there exists a canonical map $\WW(\Phi(H),\Phi(P))\to\Phi(\WW(H,P))$ of operads over $\Phi(P)$. 

This map is an isomorphism if $\Phi$ is strong symmetric monoidal. It is a weak equivalence if $P$ is a well-pointed $\Sg$-cofibrant operad, and either $H$ is an interval and $\Phi$ is a left Quillen functor, or $H$ is a cofibrant segment and $\Phi$ is an h-monoidal left Quillen functor.\end{prp}

\begin{proof}The comparison map is constructed by induction on the natural filtration of the $\WW$-construction. Since the different filtration stages $\WW_k(H,P)$ are obtained by the pushout diagrams (\ref{W}) of Section 4, and since $\Phi$ preserves colimits, the operad $\Phi(\WW(H,P))$ may be constructed as a sequential colimit of the images of (\ref{W}) under $\Phi$. For $k=0$, we have $\WW_0(\Phi(H),\Phi(P))=\Phi(P)=\Phi(\WW_0(H,P))$. For the inductive step, we consider for each $T\in\TT(n,k)$ the following commutative square of $\Aut(T)$-equivariant maps:\begin{gather}\label{h-mon}\begin{diagram}[small,silent](\Phi(H)\otimes\Phi(P))^-(T)&\rTo&\Phi((H\otimes P)^-(T))\\\dTo&&\dTo\\\Phi(H)(T)\otimes\underline{\Phi(P)}(T)&\rTo&\Phi(H(T)\otimes\underline{P}(T))\end{diagram}\end{gather}The vertical maps are given by the segment (resp. operad) structures of $\Phi(H)$ and $\,H$ (resp. $\Phi(P)$ and $\,P$), the horizontal maps are given by the symmetric monoidal structure of $\Phi$. The commutativity follows from the definitions. This allows us to define the comparison map $\WW(\Phi(H),\Phi(P))\to\Phi(\WW(H,P))$ by induction on $k$. Since the operad structure of $\WW(H,P)$ is characterised by diagram (\ref{operad}) of Section 4, the comparison map preserves the operad structure. The augmentation $\gamma_P:W(H,P)\to P$ is induced by composition maps $m_T^P:\underline{P}(T)\to P(n)$ coming from the operad structure of $P$. The analogous maps $m_T^{\Phi(P)}:\underline{\Phi(P)}(T)\to\Phi(P)(n)$ factor through $\Phi(m_T^P)$ by definition of the operad structure of $\Phi(P)$. Therefore, $\gamma_{\Phi(P)}$ factors through $\Phi(\gamma_P)$.

In the case where $\Phi$ is strong symmetric monoidal, the horizontal maps of (\ref{h-mon}) are isomorphisms and hence so are the maps $\WW_k(\Phi(H),\Phi(P))\to\Phi(\WW_k(H,P))$ as well as the colimit $\WW(\Phi(H),\Phi(P))\to\Phi(\WW(H,P))$. 

Assume now that $P$ is well-pointed $\Sg$-cofibrant, that $H$ is an interval, and that $\Phi$ a left Quillen functor. Theorem \ref{main} then implies that the augmentations $\gamma_P$ and $\gamma_{\Phi(P)}$ are weak equivalences. It follows that $\Phi(\gamma_P)$ is a weak equivalence as well, and therefore, by the 2-out-of-3 property of weak equivalences, the comparison map $\WW(\Phi(H),\Phi(P))\to\Phi(\WW(H,P))$ is a weak equivalence. 

Assume finally that $H$ is a cofibrant segment and $\Phi$ is  h-monoidal. Then the lower horizontal map of (\ref{h-mon}) is a weak equivalence. The upper horizontal map of (\ref{h-mon}) is obtained as a natural transformation of pushouts satisfying the hypothesis of Reedy's patching lemma. Therefore, the upper horizontal map is also a weak equivalence. Since in (\ref{h-mon}) all objects are $\Aut(T)$-cofibrant, and the two horizontal maps are weak equivalences, we can apply induction along $\Aut(T)\to\Sg_n$ so that at each filtration stage the  comparison map constructed above is a weak equivalence. Since the inclusions of the different filtration stages are cofibrations of cofibrant collections, we obtain a weak equivalence in the colimit by Reedy's telescope lemma.\end{proof}

At several places in this section the following categorical lemma has been used:

\begin{lma}\label{cube}Assume that the top and bottom faces of the commutative cube \begin{diagram}[width=0.5cm,height=0.5cm,silent]&&C&\rTo&&&C'\\&\ruTo&\vLine&&&\ruTo&\dTo\\A&&&\rTo&A'&&\\\dTo&&\dTo&&\dTo&&\\&&D&\hLine&\VonH&\rTo&D'\\&\ruTo&&&&\ruTo&\\B&\rTo&&&B'&&\end{diagram}are pushouts. Then, the pushout comparison maps of the left and right hand faces induce a pushout diagram\begin{diagram}[small,silent]B\cup_AC&\rTo&B'\cup_{A'}C'\\\dTo&&\dTo\\D&\rTo&D'\end{diagram}\end{lma}

\begin{proof}Use twice the fact that if a pushout square factors into two commutative squares such that the first one is a pushout, then so is the second one.\end{proof}

\section{A relative $\WW$-construction}

In this section, we will use the methods of the previous three sections
to give a \emph{relative} version of the Boardman-Vogt resolution in any monoidal model category $\EE$ with an interval $H$, as before. This
relative version applies to a suitable map of operads $P\to Q$, and
produces an interpolating operad $\WW(H,Q)_P$. The idea is that algebras
for this operad satisfy the equations for operations from $Q$ up to coherent homotopy, while they satisfy the equations for operations from $P$ on the
nose.

For a more precise statement, consider a $\Sg$-cofibration $u:P\to Q$
between well-pointed $\Sg$-cofibrant operads. Recall from \cite[Appendix]{BM}
the construction of the free extension $P[u]$ of $P$ by $u$. This free
extension is determined by the universal property that operad maps out
of $P[u]$ are in one-to-one correspondence with maps of collections out of
$Q$, whose restriction to $P$ (along $u$) is an operad map. In particular,
the identity on $Q$ corresponds to a factorisation of $u$ into operad maps
$P\to P[u]\to Q$. The relative Boardman-Vogt construction $\WW(H,Q)_P$
factors this last map $P[u]\to Q$ as expressed by the following theorem: 

\begin{thm}\label{relative}Let $\EE$ be a cofibrantly generated monoidal model category with cofibrant unit $I$ and interval $H$. Any $\Sg$-cofibration  $u:P\to Q$ of well-pointed $\Sg$-cofibrant operads factors into cofibrations $P\ito P[u]\ito\WW(H,Q)_P$ followed by a weak equivalence $\WW(H,Q)_P\eqv Q$ in such a way that the operad $\WW(H,Q)_P$ is a quotient of the operad $\WW(H,Q)$.\end{thm}

 Proposition \ref{cofibrant2} gives some evidence for Theorem \ref{relative} insofar as in the commutative square\begin{diagram}[small]\WW(H,P)&\rTo^{\WW(H,u)}&\WW(H,Q)\\\dTo^\sim&&\dTo^\sim\\P&\rTo^u&Q\end{diagram}the upper horizontal map is a cofibration of operads. Therefore, if the category of operads in $\EE$ were a \emph{left proper model category}, we could define the required factorisation simply by taking a pushout: $P\ito P\cup_{\WW(H,P)}\WW(H,Q)\eqv Q$. However, it is unclear whether the category of operads is left proper, even in the classical cases. The explicit construction of $\WW(H,Q)_P$ below is different and divides out more than the pushout just described; the quotient map $\WW(H,Q)\to\WW(H,Q)_P$ identifies \emph{any} $Q$-labelled tree with the one obtained by putting edge-lengths to $0$ whenever the two vertices of the edge are labelled by elements of $P$.\vspace{1ex}

In \cite[5.11]{BM}, the operad $P[u]$ is constructed as a sequential colimit in the category of collections; this colimit can be identified with $\WW(H,Q)_P$ in the case where $H$ is the segment $I\sqcup I$. In the absolute case where $P$ is the initial operad and $P[u]$ is the free operad on $Q$, this identification is the familiar one of Lemma \ref{counit}.\vspace{1ex}

\noindent\emph{Proof of Theorem \ref{relative}.} The operad $\WW(H,Q)_P$ is constructed as a sequential colimit of trivial cofibrations of collections:$$\WW_0(H,Q)_P\ito\WW_1(H,Q)_P\ito\WW_2(H,Q)_P\ito\cdots$$For each $k$, $\WW_k(H,Q)_P$ is a quotient of $\WW_k(H,Q)$; informally, this quotient is obtained by contracting edges with vertices on both sides labelled by elements of $P$ (no matter the length of this edge). Thus, $\WW(H,Q)_P$ is a quotient of $\WW(H,Q)$ by a filtration-preserving map.

Recall from Sections $3$ and $4$ that, for each tree $T$, we defined a trivial cofibration $H^-(T)\ito H(T)$ and an $\Aut(T)$-cofibration $\underline{Q}^-(T)\ito\underline{Q}(T)$. For the present purpose, the latter has to be replaced by another $\Aut(T)$-cofibration, written $\underline{Q}^-_P(T)\ito\underline{Q}(T)$; accordingly, the pushout $(H^-(T)\otimes \underline{Q}(T))\cup( H(T)\otimes\underline{Q}^-_P(T))$ will be denoted by $(H\otimes Q)^-_P(T)$. 

It is essential that the pushout-product map $(H\otimes Q)^-_P(T)\ito H(T)\otimes\underline{Q}(T)$ is a trivial $\Aut(T)$-cofibration. By Lemma \ref{eq2}, this will be the case if $\underline{Q}^-_P(T)\ito\underline{Q}(T)$ is an $\Aut(T)$-cofibration. Our definition of $\underline{Q}^-_P(T)$ proceeds by induction on trees and uses an intermediate object $\underline{Q}_P^r(T)$ such that, for each tree $T$, we obtain canonical $\Aut(T)$-cofibrations\begin{gather}\label{rel1}\underline{Q}^-_P(T)\ito \underline{Q}^r_P(T)\ito \underline{Q}(T).\end{gather}

Informally, $\underline{Q}^-_P(T)$ is the subobject of $\underline{Q}(T)$ formed by those labelled trees of $\underline{Q}(T)$ for which there exists either a vertex labelled by $1\in Q(1)$ or two adjacent vertices labelled by elements of $P$, while $\underline{Q}^r_P(T)$ is the subobject of those labelled trees which belong to $\underline{Q}^-_P(T)$ or have the root-vertex labelled by an element of $P$.

To start with the induction, for $T=t_1$ (resp. $T=t_n, n\not=1$), the sequence (\ref{rel1}) is defined to be $I\ito P(1)\ito Q(1)$ (resp. $0\ito P(n)\ito Q(n)$); the maps involved are $\Aut(t_n)$-cofibrations, since $P$ and $Q$ are assumed to be well-pointed and $\Sg$-cofibrant, and $u$ a $\Sg$-cofibration. 

Assume now that $T=t_n(T_1,\dots,T_n)$, and that sequences of $\Aut(T_k)$-cofibrations $\underline{Q}^-_P(T_k)\to \underline{Q}^r_P(T_k)\to \underline{Q}(T_k)$ have been defined for $k=1,\dots,n$. 

If $n=1$, we put $\underline{Q}^-_P(T)=(I\otimes\underline{Q}(T_1))\cup( P(1)\otimes\underline{Q}_P^r(T_1))\cup( Q(1)\otimes\underline{Q}_P^-(T_1))$ and $\underline{Q}^r_P(T)=(P(1)\otimes\underline{Q}(T_1))\cup (Q(1)\otimes\underline{Q}_P^-(T_1))$. By Lemma \ref{eq5}, applied to $I\ito P(1)\ito Q(1)$ and $\underline{Q}^-_P(T_1)\to \underline{Q}^r_P(T_1)\to \underline{Q}(T_1)$, we get the sequence of $\Aut(T)$-cofibrations $\underline{Q}^-_P(T)\to \underline{Q}^r_P(T)\to \underline{Q}(T)$, as required. 

Similarly, if $n\not=1$, we apply Lemma \ref{eq5} to the sequence of $\Sg_T$-cofibrations $0\ito P(n)\ito Q(n)$ (cf. Lemma \ref{eq1}) and to the sequence of $\Gamma_T$-cofibrations\begin{align*}\bigcup_{k=1}^n\underline{Q}(T_1)\otimes\cdot\otimes\underline{Q}^-_P(T_k)\otimes\cdot\otimes\underline{Q}(T_n)&\ito\bigcup_{k=1}^n\underline{Q}(T_1)\otimes\cdot\otimes\underline{Q}^r_P(T_k)\otimes\cdot\otimes\underline{Q}(T_n)\\&\ito\quad \underline{Q}(T_1)\otimes\cdots\otimes\underline{Q}(T_n)\end{align*} in order to define the sequence of $\Aut(T)$-cofibrations $\underline{Q}^-_P(T)\ito \underline{Q}^r_P(T)\ito \underline{Q}(T)$.

Now, $\WW(H,Q)_P$ is constructed by taking successive pushouts, like in the absolute case. We put $\WW_0(H,Q)_P(n)=Q(n)$ and define $\WW_k(H,Q)_P$ inductively by the following pushout (cf. Section 3 Diagram (\ref{W})):\begin{diagram}[small,silent]\coprod_{[T],T\in\TT(n,k)}(H\otimes Q)_P^-(T)\otimes_{\Aut(T)}I[\Sg_n]&\rTo^{\coprod\beta^-_T}&\WW_{k-1}(H,Q)_P(n)\\\dTo&&\dTo\\\coprod_{[T],T\in\TT(n,k)}(H(T)\otimes\underline{Q}(T))\otimes_{\Aut(T)}I[\Sg_n]&\rTo^{\coprod\beta_T}&\WW_k(H,Q)_P(n)\end{diagram}The definition of the attaching maps $\beta_T^-$ needs some care; informally, $\beta_T^-$ puts edge-lengths to $0$, whenever the vertices of the edge are both labelled by elements of $P$, and then applies the corresponding attaching map $\al^-_T$ of the absolute $\WW$-construction. That this can be done consistently follows by induction on trees, on the basis of the inductive definition of $\underline{Q}^-_P(T)$. 

The operad structure on $\WW(H,Q)_P$ is defined in a way similar to the one on $\WW(H,Q)$, cf. Section 4 Diagram (\ref{operad}). In fact, there are quotient maps $$\WW_k(H,Q)(n)\to\WW_k(H,Q)_P(n)$$ which induce a map of operads $\WW(H,Q)\to\WW(H,Q)_P$. Also, the operad structure is such that, although the trivial cofibration $Q=\WW_0(H,Q)_P\overset{\sim}{\ito}\WW(H,Q)$ is not a map of operads, its restriction to $P$ is. Therefore, we get an operad map $P[u]\to\WW(H,Q)_P$ for which the following diagram commutes:\begin{diagram}[small,silent]P[u]&\rTo&\WW(H,Q)_P\\\uTo&&\dTo^\sim\\P&\rTo^u&Q\end{diagram}The right vertical map is induced by the counit $\eps:H\to I$ and the operad structure maps $\underline{Q}(T)\to Q(n)$; it is a weak equivalence since the trivial cofibration $\WW_0(H,Q)_P\overset{\sim}{\ito}\WW(H,Q)_P$ induces a section on the underlying collections. Since the free extension $P\to P[u]$ is a cofibration of operads, it remains to be shown that the upper horizontal map $P[u]\to\WW(H,Q)_P$ is a cofibration of operads. For this, we use that$$P[u]=\WW(I\sqcup I,Q)_P$$can be constructed similarly as a sequential colimit$$Q=P[u]_0\ito P[u]_1\ito P[u]_2\ito \dots$$and we refine the filtration of $\WW(H,Q)_P$ exactly as in the absolute case:\begin{diagram}[small,silent]\WW_0(H,Q)_P&\rTo&\WW_1(H,Q)_P&\rTo&\WW_2(H,Q)_P&\rTo&\WW_3(H,Q)_P\,\cdots\\\dTo&\ruTo&\dTo&\ruTo&\dTo&\ruTo\\\WW_0^+(H,Q)_P&&\WW_1^+(H,Q)_P&&\WW^+_2(H,Q)_P&\cdots\end{diagram}Then $P[u]\to\WW(H,Q)_P$ restricts to $P[u]_{k+1}\ito\WW_k^+(H,Q)_P$ on the different filtration layers, and we can argue like in the proof of Lemma \ref{khom}.\qed

\section{Particular instances of the $\WW$-construction}

\subsection{The Boardman-Vogt resolution of topological operads}

The set-theore-tical description of Section $1$ applies to topological operads and is clearly equivalent to Boardman and Vogt's original definition \cite[Chapter III]{BV} of the $\WW$-construction, as well as to our inductive definition of Section $4$. It is interesting to observe the different choices of segment structures in the literature. Boardman and Vogt endow the unit-interval $[0,1]$ with ``reversed'' multiplication; i.e. the product $s\star t$ is defined by the identity $(1-s\star t)=(1-s)(1-t)$ which makes $0$ neutral and $1$ absorbing; they already observe that other choices are possible, cf. \cite[Remark 3.2]{BV}. In Vogt's recent work \cite{V}, the unit-interval $[0,1]$ is endowed with the maximum-operation, which is the natural choice from our point of view. Kontsevich and Soibelman \cite{KS} choose $[0,\infty]$ as segment: addition makes $0$ neutral and $\infty$ absorbing.

The homotopy invariance of $\WW([0,1],P)$-structures along homotopy equivalences is proved directly by Boardman and Vogt \cite[Theorem 4.37]{BV}, and has probably been one of their main reasons for introducing the $\WW$-construction. Using Theorem \ref{main}, the homotopy invariance of $\WW([0,1],P)$-structures follows from \cite[Theorem 3.5]{BM} by the cofibrancy of $\WW([0,1],P)$, whenever the topological operad $P$ is well-pointed $\Sg$-cofibrant. Observe that there are two monoidal model structures on compactly generated spaces: Quillen's structure \cite{Q} with weak homotopy equivalences as weak equivalences, and Str\o m's structure \cite{S} with genuine homotopy equivalences as weak equivalences. However, only Quillen's structure is cofibrantly generated, a hypothesis we made in the preliminary Section 2.5 on equivariant monoidal homotopy theory. Therefore, some care is needed to adapt our methods to the non-cofibrantly generated case. The cofibrancy of $\WW([0,1],P)$ with respect to Str\o m's model structure is precisely Vogt's result \cite[Theorem 4.1]{V}.

\subsection{The Boardman-Vogt resolution of simplicial operads}\label{simplicial}

The category of simplicial sets is one of the most basic cofibrantly generated monoidal model categories, cf. the Appendix for our notation and terminology. Its cofibrations are the monomorphisms, its weak equivalences are the realisation weak equivalences, and its fibrations are the Kan fibrations, cf. Quillen \cite{Q}. The unit of $\Sets^\Dop$ is the representable $0$-simplex $\Delta^0=\Delta(-,[0])$. The representable $1$-simplex $\Delta^1=\Delta(-,[1])$ is an interval; indeed, the segment structure $\Delta^0\sqcup\Delta^0\to\Delta^1\to\Delta^0$ is induced by simplicial coface and codegeneracy operators; the binary operation $\Delta^1\times\Delta^1\to\Delta^1$ is induced by the maximum-operation $[1]\times[1]\to[1]$. Therefore, Theorem \ref{main} implies that each well-pointed $\Sg$-cofibrant simplicial operad $P$ has a cofibrant resolution $\WW(\Delta^1,P)$. All simplicial operads are well-pointed. The $\Sg$-cofibrant simplicial operads are precisely those with free $\Sg_n$-actions.

The geometric realisation functor $|\!-\!|:\Sets^\Dop\to\Top$ is a strong symmetric monoidal left Quillen functor. Proposition \ref{Quillen} thus implies that the realisation $|\WW(\Delta^1,P)|$ of the simplicial $\WW$-construction is isomorphic to the topological $\WW$-construction $\WW(|\Delta^1|,|P|)$. Observe that $|\Delta^1|$ is the unit-interval $[0,1]$ endowed with the maximum-operation.

An important example of a $\Sg$-cofibrant simplicial operad is the simplicial $E_\infty$-operad $E\Sg$ formed by the homogeneous bar resolutions $E\Sg_n$ of the symmetric groups $\Sg_n,\,n\geq 0$. This $E_\infty$-operad admits a canonical filtration by simplicial $E_k$-suboperads $E_k\Sg$, cf. \cite{B}. The simplicial operad $\WW(\Delta^1,E\Sg)$ is thus a \emph{cofibrant $E_\infty$-operad} which comes equipped with a filtration by \emph{cofibrant $E_k$-suboperads} $\WW(\Delta^1,E_k\Sg)$. Moreover, Proposition \ref{cofibrant2} shows that the inclusions of the different layers $\WW(\Delta^1,E_k\Sg)\to\WW(\Delta^1,E_{k+1}\Sg)$ are cofibrations of operads.

\subsection{The Godement resolution of an operad}\label{Godement0}

The adjunction between pointed collections and operads in a symmetric monoidal category $\EE$ induces for each operad $P$ in $\EE$ a well known simplicial resolution$$\Delta^\op\to\Oper(\EE):[k]\mapsto\GG_k(P)$$where $\GG_k(P)=(\FF_*\UU_*)^{k+1}(P)$, and where the simplicial operators are generated by \begin{align*}\partial_i=(\FF_*\UU_*)^{k-i}\eps_{(\FF_*\UU_*)^iP}&:\GG_k(P)\to\GG_{k-1}(P)\\s_i=(\FF_*\UU_*)^{k-i}\FF_*\eta_{\UU_*(\FF_*\UU_*)^iP}&:\GG_k(P)\to\GG_{k+1}(P)\end{align*}for $i=0,\dots,k$. Here $\eps:\FF_*\UU_*\Rightarrow id_{\Oper(\EE)}$ (resp. $\eta:id_{\Coll_*(\EE)}\Rightarrow\UU_*\FF_*$) is the counit (resp. unit) of the adjunction $\FF_*:\Coll_*(\EE)\lrto\Oper(\EE):\UU_*$. The simplicial operad $\GG_\cdot(P)$ is augmented over $P$ by the counit $\eps_P:\GG_0(P)\to P$. After application of the forgetful functor, we get a simplicial object in pointed collections which contains the underlying collection of $P$ as a simplicial deformation retract. This justifies the terminology ``simplicial resolution''. In the literature, the latter is often called ``cotriple resolution'' or ``comonadic resolution''. It is a special case of May's \emph{double-sided bar} construction \cite{M}. For the purpose of the present article, we refer to $\GG_\cdot(P)$ as the \emph{Godement resolution} of $P$ since, to our knowledge, the simplicial formulae above occur for the first time in Godement's Appendix of \cite{G}.

The unit $I$ of $\EE$ defines a functor $(-)_\EE:\Sets\to\EE:S\mapsto\coprod_SI$. The latter is strong symmetric monoidal since the monoidal structure of $\EE$ distributes over coproducts by closedness. Post-composition with $(-)_\EE$ yields a strong symmetric monoidal functor $\Sets^\Dop\to\EE^\Dop$, also denoted by $(-)_\EE$. To be explicit, for a simplicial set $X$, we have $(X_\EE)_n=(X_n)_\EE$. Here, we use that the (closed) symmetric monoidal structure of $\EE^\Dop$ is defined degreewise. Moreover, if $\EE$ is a monoidal model category, then so is $\EE^\Dop$ with respect to the \emph{Reedy model structure} on simplicial objects, cf. the Appendix. The functor $(-)_\EE$ takes cofibrations of simplicial sets to cofibrations of simplicial objects in $\EE$. 

It follows that the category $\EE^\Dop$ has a canonical cofibrant segment, namely the image $\Delta^1_\EE$ of the standard simplical interval $\Delta^1$, cf. Section \ref{simplicial}. Notice however that $\Delta_\EE^1$ is only a cofibrant segment, but not an interval for the Reedy model structure on $\EE^\Dop$, since the counit $\Delta^1_\EE\to\Delta^0_\EE$ is not a weak equivalence with respect to the Reedy model structure. There is a canonical way of rendering an arbitrary map of a cofibrantly generated model category into a weak equivalence (thus enlarging the class of weak equivalences) while keeping fixed the class of cofibrations; this process is called a \emph{left localisation} of the model category. We refer the reader to the book of Hirschhorn \cite{Hir} for more details. We consider $\EE$ as the full subcategory of $\EE^\Dop$ formed by the constant simplicial objects.

\begin{sprp}\label{Godement}For each operad $P$ in $\EE$, there is a canonical isomorphism$$\GG_\cdot(P)\cong\WW(\Delta^1_\EE,P)$$of $P$-augmented simplicial operads in $\EE$. In particular, the Godement resolution $\GG_\cdot(P)$ is a cofibrant simplicial resolution for well-pointed $\Sg$-cofibrant operads $P$ in any left localisation of $\EE^\Dop$ in which the counit $\Delta^1_\EE\to\Delta^0_\EE$ is a weak equivalence.\end{sprp}

\begin{proof}Applying inductively Lemma \ref{diamond}, we obtain $$\GG_k(P)\cong\WW(I^{\diamond (k+1)},P)=\WW(I^{\sqcup (k+2)},P).$$ The face operator $\partial_i:\GG_k(P)\to\GG_{k-1}(P)$ is induced by the codiagonal $I\sqcup I\to I$ applied to the $(i+1)$th and $(i+2)$th summand of $I^{\sqcup (k+2)}$. Indeed, Lemma \ref{diamond} gives this identification for the last face operator $\partial_k$; for the other face operators, the identification follows by induction on $k$. A similar argument shows that the degeneracy operator $s_i:\GG_k(P)\to\GG_{k+1}(P)$ is induced by $0:I\to I\sqcup I$ applied to the $(i+1)$th summand of $I^{\sqcup (k+2)}$.

On the other hand, the object of $k$-simplices of $\Delta_\EE^1$ is given by $$(\Delta_\EE^1)_k=\Delta([k],[1])_\EE=I^{\sqcup(k+2)}$$with the same simplicial operators as above and with segment structure$$(\Delta_\EE^0)_k\sqcup(\Delta_\EE^0)_k\to(\Delta_\EE^1)_k\to(\Delta_\EE^0)_k$$identical to the canonical one on $I^{\diamond (k+1)}=I^{\sqcup(k+2)}$. We thus get isomorphisms $\GG_k(P)\cong\WW((\Delta^1_\EE)_k,P)$ compatible with the simplicial structures on both sides. 
The second assertion follows from Theorem \ref{main}, since the embedding of $\EE$ in $\EE^\Dop$ preserves well-pointed $\Sg$-cofibrant operads, and since $\Delta^1_\EE$ is an interval in any left localisation of $\EE^\Dop$ in which $\Delta_\EE^1\to\Delta_\EE^0$ is a weak equivalence.\end{proof}

\subsection{Comparison of the Boardman-Vogt and Godement resolutions}\label{comparison}

In order to compare the Boardman-Vogt and Godement resolutions of an operad, the latter has to be realised in the ambient monoidal model category $\EE$. It is shown in the Appendix that a convenient way of doing so, is to realise the Godement resolution with respect to a so-called \emph{standard system of simplices} in $\EE$. This is a \emph{cosimplicial object} $C:\Delta\to\EE$ which is \emph{h-monoidal, Reedy-cofibrant} and \emph{h-constant}; see the Appendix for details. The $1$-truncation of a standard system of simplices $C$ defines an interval $C^0\sqcup C^0\to C^1\to C^0$, cf. Lemma \ref{interval2}. It is therefore natural to compare $\WW(C^1,P)$ and $|\GG_\cdot(P)|_{C}$. We give an explicit comparison map from the former to the latter. Geometrically, this comparison map is a kind of Eilenberg-Zilber map which ``subdivides'' the cubical Boardman-Vogt resolution so as to get the realisation of the simplicial Godement resolution. 

In the topological case, the existence of such a comparison map has been known for a while; in particular, the coexistence of a cubical and a simplicial description of the \emph{same} topological object is a recurring theme in the theory of homotopy colimits, cf. especially the work of Cordier and Porter \cite{CP}.

\begin{sthm}For any operad $P$ of a cofibrantly generated monoidal model category $\EE$ with monoidal cosimplicial object $C$, there is a canonical comparison map $\WW(C^1,P)\to|\GG_\cdot(P)|_{C}$ of $P$-augmented operads in $\EE$. 

This comparison map is an isomorphism if $C$ is strong monoidal; it is a weak equivalence if $P$ is well-pointed $\Sg$-cofibrant and $C$ is h-monoidal Reedy-cofibrant; it is a weak equivalence of resolutions if $P$ is well-pointed $\Sg$-cofibrant and $C$ is a standard system of simplices.\end{sthm}

\begin{proof}The simplicial segment $\Delta^0_\EE\sqcup\Delta^0_\EE\to\Delta^1_\EE\to\Delta_\EE^0$ is cofibrant. Its image under the realisation functor $|\!-\!|_C$ is the segment $C^0\sqcup C^0\to C^1\to C^0$ derived from the given monoidal cosimplicial object $C$, cf. Lemma \ref{bi}. Moreover, the functor $|\!-\!|_{C}:\EE^\Dop\to\EE$ is strong symmetric monoidal (resp. an h-monoidal left Quillen functor) whenever the cosimplicial object $C$ is strong monoidal (resp. h-monoidal Reedy-cofibrant). Therefore, Propositions \ref{Quillen} and \ref{Godement} give a comparison map$$\WW(C^1,P)\cong\WW(|\Delta_\EE^1|_{C},|P|_{C})\to|\WW(\Delta_\EE^1,P)|_{C}\cong|\GG_\cdot(P)|_{C}$$ which is an isomorphism in the strong monoidal case, and a weak equivalence in the h-monoidal Reedy-cofibrant case. If $C$ is a standard system of simplices, its $1$-truncation is an interval by Lemma \ref{interval2}, whence $\WW(C^1,P)$ (by Theorem \ref{main}) and $|\GG_\cdot(P)|_{C}$ (by the 2-out-of-3-property of weak equivalences) are resolutions of $P$.\end{proof}

\begin{srmk}The preceding theorem applies in particular to simplicial, topological, spectral and chain operads, cf. Examples \ref{standard}. In the first three cases we get an \emph{isomorphism} between the Boardman-Vogt resolution and the realisation of the Godement resolution; in the fourth case we obtain a \emph{weak equivalence} of resolutions.\end{srmk}

\subsection{The cobar-bar resolution of reduced chain operads}

The category of $\ZZ$-graded chain complexes $\Ch(R)$ over a unitary commutative ring $R$ is a monoidal model category for the projective model structure, i.e., the fibrations are precisely the epimorphisms, and the weak equivalences are precisely the quasi-isomorphisms, cf. Hovey \cite{Hov}. This monoidal model category has an interval, namely the image $\NN^R_*(\Delta^1)$ of the simplicial interval $\Delta^1$ under the normalised $R$-chain functor $\NN^R_*$. This implies in particular that the category of reduced $R$-chain operads carries a canonical model structure, cf. \cite[3.3.3]{BM}. The purpose of this section is to show that for any reduced $R$-chain operad $P$, the reduced $\WW$-construction $\WW^{red}(\NN^R_*(\Delta^1),P)$ is \emph{isomorphic} to the cobar-bar resolution $\BB^c(\BB(P))$ of $P$. By Theorem \ref{main}, this recovers Theorem 3.2.16 of Ginzburg and Kapranov \cite{GK}, with the refinement that the cobar-bar resolution of $P$ is shown to be \emph{cofibrant} whenever $P$ is $\Sg$-cofibrant.

\subsubsection{The canonical $R$-chain interval}

The normalised $R$-chain functor $\NN^R_*$ is an h-monoidal left Quillen functor from simplicial sets to $R$-chain complexes, cf. Remark \ref{standard}. Therefore, $\NN_*^R(\Delta^1)$ is an interval in $\Ch(R)$ by Corollary \ref{hfc} and Lemma \ref{interval2}. Explicitly, in each degree $\NN_*^R(\Delta^1)$ is freely generated by the set of non-degenerate simplices of $(\Delta^1)_*$; therefore, there are two generators $\gamma_0,\gamma_1$ for $\NN_0^R(\Delta^1)$, and one generator $\gamma$ for $\NN_1^R(\Delta^1)$, which are related by the differential $$\partial(\gamma)=\gamma_1-\gamma_0.$$The binary operation $\vee:\NN^R_*(\Delta^1)\otimes_R\NN^R_*(\Delta^1)\to\NN_*^R(\Delta^1)$ is determined by the requirement that $\gamma_0$ is neutral, $\gamma_1$ is absorbing; in particular, for degree reasons we have the identities $\gamma\vee\gamma_1=\gamma_1\vee\gamma=0$ and $\gamma\vee\gamma=0$.

\subsubsection{The $\WW$-construction for reduced $R$-chain operads}We shall give an explicit description of $\WW^{red}(H,P)$ for $H=\NN^R_*(\Delta^1)$ and $P$ a reduced $R$-chain operad. All tensors are to be taken over $R$. The abelian structure of $\Ch(R)$ allows us to express $\WW^{red}(H,P)$ in terms of the corresponding $\WW$-construction for pseudo-operads, cf. Remark \ref{analog}(b) and (c). Indeed, the unit of any reduced operad $P$ has a retraction by $P(1)\cong P(1)\otimes P(0)\to P(0)\cong I.$ Therefore, the unit of a reduced $R$-chain operad splits off as $P(1)=R\oplus\bar{P}(1)$. We extend $\bar{P}(1)$ to a pseudo-operad, by setting $\bar{P}(0)=0$ and $\bar{P}(n)=P(n)$ for $n\geq 2$. In fact, $\bar{P}$ is the unique pseudo-sub-operad of $P$ such that $\underline{R}\oplus\bar{P}=P$, where $\underline{R}(0)=\underline{R}(1)=R$ and $\underline{R}(n)=0$ for $n\geq 2$. The reduced $\WW$-construction in $\Ch(R)$ simplifies then to $$\WW^{red}(H,P)=\underline{R}\oplus\WW^{ps}(H,\bar{P}).$$where $\WW^{ps}(H,\bar{P})$ is defined as in Section $4$ except that the unit-tree and all trees with vertices of valence $0$ are discarded, and no identifications for units are made. In particular, everywhere the map $(H\otimes\underline{P})^-(T)\to H(T)\otimes\underline{P}(T)$ is replaced by the map $H^-(T)\otimes\underline{\bar{P}}(T)\to H(T)\otimes\underline{\bar{P}}(T).$ As a consequence we can also ignore the multiplicative structure of $H$.

For an $R$-chain complex $(A,\partial^A)$, we define a shifted $R$-chain complex $(A[d],\partial^{A[d]})$ by $A[d]_i=A_{i-d}$ and $\partial_i^{A[d]}=(-1)^d\partial_{i-d}^A.$ In particular, the Koszul rule for signs gives a canonical identification $A[d]= R[d]\otimes A$. For each tree $T$ with set $E(T)$ of internal edges, we then have (ignoring differentials)\begin{align*}H(T)&=\bigoplus_{E(T)=E_{\gamma_0}\sqcup E_{\gamma_1}\sqcup E_\gamma} R^{\otimes|E_{\gamma_0}|}\otimes R^{\otimes|E_{\gamma_1}|}\otimes R[1]^{\otimes|E_\gamma|}\\&= \bigoplus_{E(T)=E_{\gamma_0}\sqcup E_{\gamma_1}\sqcup E_\gamma}R[|E_\gamma|].\end{align*}Colimits preserve sums, so we can write each $\WW_k^{ps}(H,\bar{P})$ as a sum. The summands corresponding to edge-decompositions with $E_{\gamma_0}\not=\emptyset$ are identified with summands in $\WW_{k-1}^{ps}(H,\bar{P})$. Therefore, we get\begin{gather}\label{cb1}\WW^{ps}(H,\bar{P})=\bigoplus_{T\not=\,|}(\bigoplus_{E_\gamma\subseteq E(T)}\underline{\bar{P}}(T)[|E_\gamma|])\otimes_{\Aut(T)}I[\Sg_{|T|}]\end{gather}as a collection of graded $R$-modules. There is a more appealing way to write the same formula. Each summand is indexed by a pair $(T,E_\gamma)$ with $E_\gamma\subseteq E(T)$. Such a pair may be considered as a ``tree of trees'': cut internal edges of $T$ not belonging to $E_\gamma$ in the middle, and consider the connected components as the vertices of a quotient tree $T/E_\gamma$; by construction, the valence of each vertex of $T/E_\gamma$ equals the number of input edges of the corresponding subtree of $T$. Therefore, for a given tree $S$, the partial sum in (\ref{cb1}) over all pairs $(T,E_\gamma)$ with $T/E_\gamma\cong S$, is isomorphic (under rearrangement of the factors) to $\underline{\FF(\bar{P}[1])}(S)[-1]$. The degree-shifts account for the fact that a tree with $d$ internal edges has $d+1$ vertices. Formula (\ref{cb1}) may thus be rewritten as follows:\begin{gather}\label{cb2}\WW^{ps}(H,\bar{P})=\bar{\FF}(\FF(\bar{P}[1])[-1]),\end{gather}where $\bar{\FF}$ is the free operad construction without unit. 

The augmentation $\WW^{red}(H,P)\to P$ is induced by the counit $\eps_H:H\to R$ which is given by $\eps_H(\gamma_0)=\eps_H(\gamma_1)=1$ and $\eps_H(\gamma)=0$: so it takes summands indexed by $(T,E_\gamma)$ to $0$ whenever $E_\gamma\not=\emptyset$, and coincides otherwise with the counit $\FF_*(P)\to P$.

Let us describe the differential $\partial^\WW$ of $\WW^{red}(H,P)$. It has two components $$\partial^\WW=\partial^P+\partial^H$$ induced respectively by the differential of $P$ and the differential of $H$. We shall call $\partial^P$ the \emph{internal} differential and $\partial^H$ the \emph{external} differential of $\WW^{red}(H,P)$. The internal differential is induced by the usual one for tensor-products applied to$$\underline{P}(T)=\bigotimes_{v\in T}P(|v|).$$Its explicit formula depends on a convention for the ordering of the vertices of $T$. We choose vertex-orderings for pairs $(T,E_\gamma)$ in such a way that the above described identification of $(T,E_\gamma)$ with a quotient-tree $T/E_\gamma$ of subtrees of $T$ is compatible with these orderings. The internal differential $\partial^P$ is then the one obtained by formula (\ref{cb2}) applied to the $R$-chain operad $P$. 

Accordingly, the external differential $\partial^H$ may be described as follows. Let the pair $(T,E_\gamma)$ correspond to a tree $T$ with $r$ subtrees $S_1,\dots,S_r$ having $\sg_1,\dots,\sg_r$ internal edges respectively. In particular, $|E_\gamma|=\sg_1+\cdots+\sg_r$. The restriction of $\partial^H$ to the summand $\underline{\bar{P}}(T)[|E_\gamma|]=\bigotimes_{i=1}^r\underline{\bar{P}}(S_i)[\sg_i]$ is the alternating sum$$\sum_{i=1}^r(-1)^{i-1}1^{\otimes(i-1)}\otimes\partial_i\otimes 1^{\otimes(r-i-1)}$$where $\partial_i$ acts on $\underline{\bar{P}}(S_i)[\sg_i]$ only. We have$$\underline{\bar{P}}(S_i)[\sg_i]=(\bigotimes_{v\in S_i}\bar{P}(|v|)[1])[-1]\quad\text{and}\quad\partial_i=\sum_{e^i_j\in E(S_i)}\partial^0_{e_j^i}-\partial^1_{e_j^i}.$$The operator $\partial_{e_j^i}^0$ (resp. $\partial_{e_j^i}^1$) gives $e_j^i$ length $\gamma_1$ (resp. length $\gamma_0$). 

More precisely, $\partial^0_{e^i_j}$ maps $\underline{\bar{P}}(S_i)[\sg_i]$ to $\underline{\bar{P}}(S_i^-)[\sg_i^-]\otimes\underline{\bar{P}}(S_i^+)[\sg_i^+]$ (where $S_i^-$ and $S_i^+$ are obtained from $S_i$ by cutting $e_j^i$ in the middle and $\sg_i^\pm$ is the number of internal edges in $S_i^\pm$ so that $\sg_i-1=\sg_i^-+\sg_i^+$) via a shift of the canonical isomorphism $\underline{\bar{P}}(S_i)\cong\underline{\bar{P}}(S_i^-)\otimes\underline{\bar{P}}(S_i^+)$; the operator $\partial^1_{e^i_j}$ maps $\underline{\bar{P}}(S_i)[\sg_i]$ to $\underline{\bar{P}}(S_i/\{e^i_j\})[\sg_i-1]$ via a shift of the operad composition map $\underline{\bar{P}}(S_i)\to\underline{\bar{P}}(S_i/\{e^i_j\})$.\vspace{1ex}

\subsubsection{The cobar-bar adjunction}In order to state our comparison theorem explicitly, we briefly recall the cobar-bar adjunction, see Getzler and Jones \cite[Section 2.1]{GJ} and Fresse \cite[Part 3]{F} for details. Up to a two-fold dualisation (which avoids cooperads), the counit of cobar-bar adjunction also coincides with the $\mathbb{D}\mathbb{D}$-resolution of Ginzburg and Kapranov \cite{GK}, cf. Markl, Shnider and Stasheff \cite[Section 3.1]{MSS}.

A reduced operad $P$ such that $P(n)=0$ for $n\geq 2$, is the same as an augmented $R$-chain algebra; dually, a reduced cooperad $C$ such that $C(n)=0$ for $n\geq 2$, is the same as a coaugmented $R$-chain coalgebra. The cobar-bar adjunction between reduced cooperads and reduced operads generalises the classical cobar-bar adjunction between coaugmented coalgebras and augmented algebras. We follow Getzler and Jones \cite{GJ} in using twisting cochains for the definition of the adjunction.

We assume throughout that $C$ is a reduced $R$-chain cooperad and $P$ is a reduced $R$-chain operad. For the existence of the adjunction no further constraint is required; however, the unit of the adjunction will not be a quasi-isomorphism unless some connectivity assumption on $C$ is made. We denote by $\bar{C}$ the unique pseudo-sub-cooperad of $C$ such that $\underline{R}\oplus\bar{C}=C$. 

A \emph{twisting cochain} $C\to P$ is a degree $-1$ map $\tau:\bar{C}\to\bar{P}$ of collections of graded $R$-modules such that $D\tau=\tau\cup\tau$. By definition $D\tau=\partial^{\bar{P}}\tau+\tau\partial^{\bar{C}}$, while the cup square $\tau\cup\tau$ is the following composite map of degree $-2$: \begin{gather}\label{twist}\bar{C}\lra\FF^{(2)}(\bar{C})\overset{\FF^{(2)}(\tau)}{\lra}\FF^{(2)}(\bar{P})\lra\bar{P}.\end{gather}Here $\FF^{(2)}$ denotes the ``quadratic'' part of the free operad functor, i.e. the summand of $\FF$ indexed by trees with exactly two vertices or, equivalently, one internal edge, cf. Remark \ref{filter}. Since the expansions of the cofree cooperad and the free operad of a collection are the same, the first map in (\ref{twist}) is induced by the cooperad structure of $C$ while the third map in (\ref{twist}) is induced by the operad structure of $P$. The set $Twist(C,P)$ of twisting cochains $C\to P$ is a contravariant (resp. covariant) functor in maps of reduced cooperads $C'\to C$ (resp. reduced operads $P\to P'$).


There is a unique coderivation $\partial^{bar}$ of the cofree cooperad $(\FF^c(\bar{P}[1]),\partial^c)$ such that the counit $\eps_{\bar{P}[1]}$ of the forgetful-cofree adjunction becomes a twisting cochain $(\FF^c(\bar{\PP}[1]),\partial^c+\partial^{bar})\to P$. Dually, there is a unique derivation $\partial^{cobar}$ of the free operad $(\FF(\bar{C}[-1]),\partial)$ such that the unit $\eta_{\bar{C}[-1]}$ of the free-forgetful adjunction becomes a twisting cochain $C\to(\FF(\bar{C}[-1]),\partial+\partial^{cobar})$. Indeed, it follows from the definitions that, up to a sign, the coderivation $\partial^{bar}$ uniquely extends the cup square $\eps_{\bar{P}[1]}\cup\eps_{\bar{P}[1]}:\FF^c(\bar{P}[1])\to\bar{P}[1]$, while the derivation $\partial^{cobar}$ uniquely extends the cup square $\eta_{\bar{C}[-1]}\cup\eta_{\bar{C}[-1]}:\bar{C}[-1]\to\FF(\bar{C}[-1])$. Since $\eps_{\bar{P}[1]}$ and $\eta_{\bar{C}[-1]}$ are non-zero only on summands indexed by trees with exactly one vertex, their cup squares are non-zero only on the quadratic parts $\FF^{(2)}(\bar{P}[1])$ resp. $\FF^{(2)}(\bar{C}[-1])$, cf. Getzler and Jones \cite[2.2-4]{GJ} and Fresse \cite[3.1.9-10]{F}. 

By definition, the bar cooperad of $P$ is given by $\BB(P)=(\FF^c(\bar{P}[1]),\partial^c+\partial^{bar})$ and the cobar operad of $C$ is given by $\BB^c(C)=(\FF(\bar{C}[-1]),\partial+\partial^{cobar})$. It turns out that $\partial^{bar}$ and $\partial^{cobar}$ are actually differentials, the so-called \emph{bar} and \emph{cobar} differentials. The cobar-bar adjunction arises then from the following commutative diagram (where $\#$ forgets differentials):\begin{diagram}[small,silent]\Oper(\FF(\bar{C}[-1])^\#,P^\#)&\llra&\Coll_{*}(C^\#,P^\#)[-1]&\llra&\Cooper(C^\#,\FF^c(\bar{P}[1])^\#)\\\uInto&&\uInto&&\uInto\\\Oper(\BB^c(C),P)&\llra&Twist(C,P)&\llra&\Cooper(C,\BB(P))\end{diagram}In particular, the counit of cobar-bar adjunction $\BB^c(\BB(P))\to P$ is induced by the above-mentionned twisting cochain $\BB(P)\to P$.

At several places in the literature it has been suggested that the cobar-bar chain resolution of Ginzburg and Kapranov \cite{GK} is analogous to the topological $\WW$-construction of Boardman and Vogt \cite{BV}, cf. Kontsevich and Soibelman \cite{KS} and also Markl, Shnider and Stasheff \cite[pg. 128f.]{MSS}. For our general $\WW$-construction, this analogy becomes an isomorphism:\begin{sthm}For each reduced $R$-chain operad $P$, there is a canonical isomorphism of $P$-augmented operads$$\WW^{red}(\NN_*^R(\Delta^1),P)\cong\BB^c(\BB(P)).$$In particular, the counit of the cobar-bar adjunction is a cofibrant resolution of $P$ whenever the collection underlying $P$ is cofibrant.\end{sthm}

\begin{proof}The second assertion follows from the first by Theorem \ref{main} and the fact that $\Sg$-cofibrant reduced operads are well-pointed. In (\ref{cb2}) we showed that the reduced $\WW$-construction $\WW^{red}(H,P)$ for $H=\NN_*^R(\Delta^1)$ may be expanded as $\FF(\FF(\bar{P}[1])[-1])$, which is precisely the expansion of $\BB^c(\BB(P))$. The differential of $\WW^{red}(H,P)$ is a sum  of two differentials $\partial^P+\partial^H$. The internal differential $\partial^P$ coincides with the ``free-cofree'' part of the differential of $\BB^c(\BB(P))$. The external differential $\partial^H$ is itself a sum of two differentials corresponding to the two vertices of $\Delta^1$. The $\partial^0$-differential coincides (up to a sign) with the part of the differential of $\BB^c(\BB(P))$ which is induced by the bar differential of $\BB(P)$, while the $\partial^1$-differential coincides (up to a sign) with the cobar differential of $\BB^c(\BB(P))$. In particular, the almost cofree cooperad $\BB(P)$ may be identified (up to a degree-shift) with the sum of those summands of (\ref{cb1}) which are indexed by pairs $(T,E_\gamma)$ such that $E_\gamma=E(T)$.\end{proof}

\section{Appendix: Geometric realisation in monoidal model categories}

In this Appendix, we collect some results on the realisation of simplicial objects in monoidal model categories. These results have been used for the comparison of the Boardman-Vogt and Godement resolutions of an operad, cf. Section \ref{comparison}. They may be useful in other contexts as well.

Throughout, $\EE$ denotes a closed symmetric monoidal category which is cocomplete and finitely complete. In most cases, $\EE$ will be a cofibrantly generated monoidal model category with cofibrant unit $I$, cf. Section 2.

The category of non-empty finite ordinals $[n],\,n\geq 0,$ and order-preserving maps will be written $\Delta$. A \emph{simplicial object} in $\EE$ is a contravariant functor $X:\Delta^{op}\to \EE$. As usual, we write $X_n$ for $X([n])$. The category $\EE^\Dop$ of simplicial objects in $\EE$ is again symmetric monoidal for the degreewise tensor product. A cosimplicial object is a covariant functor $C:\Delta\to\EE$. We shall write $C^n$ for $C([n])$. Each cosimplicial object $C$ defines a \emph{realisation functor} $|\!-\!|_{C}:\EE^\Dop\to\EE$, given on objects by the coend formula $|X|_{C}=X_\cdot\otimes_\Delta C^\cdot$. Explicitly, this coend is the coequaliser of the following pair in $\EE$:$$\coprod_{\phi:[m]\to[n]}(\phi^*\otimes id_{C^m},id_{X_n}\otimes\phi_*):\coprod_{\phi:[m]\to[n]}X_n\otimes C^m\dto\coprod_{[n]}X_n\otimes C^n$$

Each set $S$ defines an object $S_\EE=\coprod_SI$ of $\EE$. This functor is strong symmetric monoidal and admits a strong symmetric monoidal prolongation to simplicial objects, also denoted$$(-)_\EE:\Sets^\Dop\to\EE^\Dop.$$Here, we use that the tensor of $\EE$ distributes over coproducts in $\EE$ by closedness. We shall use the symbol $\Delta^n$ for the representable $n$-simplex in $\Sets^\Dop$. Accordingly, its image in $\EE^\Dop$ will be denoted by $\Delta_\EE^n$. Since $(-)_\EE$ is strong monoidal, we have canonical isomorphisms $\Delta^m_\EE\otimes\Delta^n_\EE\cong (\Delta^m\times\Delta^n)_\EE$ for $m,n\geq 0.$


\begin{dfn}\label{defmon}A cosimplicial object $C$ in $\EE$ is \emph{(strong) monoidal}  if the realisation functor $|\!-\!|_C:\EE^\Dop\to\EE$ is unit-preserving and (strong) symmetric monoidal.\end{dfn}

Our first objective is to show that these apparently global properties of a cosimplicial object of $\EE$ are already determined by the behaviour of the realisation functor on tensors of the form $\Delta^m_\EE\otimes\Delta^n_\EE$. (The argument does not actually involve any particular feature of the category $\Delta$, and an analogous property holds for $\EE$-valued presheaves over any small category).

For any two simplicial objects $X,Y$ of $\EE$, we write $X\Box Y$ for the bisimplicial object $(\Delta\times\Delta)^\op\to\EE:([m],[n])\mapsto X_m\otimes Y_n$. The diagonal of $X\Box Y$ is $X\otimes Y$.

Similarly, for cosimplicial objects $C,D$ of $\EE'$, we write $C\Box D$ for the bicosimplicial object $\Delta\times\Delta\to\EE':([m],[n])\mapsto C^m\otimes D^n$. In particular, there is such an object $\Delta_\EE\Box\Delta_\EE:\Delta\times\Delta\to\EE^\Dop$ for the cosimplicial object $\Delta_\EE:\Delta\to\EE^\Dop:[n]\mapsto\Delta^n_\EE$.

Finally, for any functor $F$ with values in $\EE^\Dop$, we write $|F|_C$ for the composite functor $|\!-\!|_C\circ F$ which takes values in $\EE$.

\begin{lma}\label{bi}For any simplicial objects $X,Y$ of $\,\EE$, and any cosimplicial object $C$ of $\EE$, there are canonical $\EE$-isomorphisms natural in $X,Y,C:$\begin{align*}|\Delta_\EE|_C&\cong C,\\(X\Box  Y)\otimes_{\Delta\times\Delta}(|\Delta_\EE|_C\Box |\Delta_\EE|_C)&\cong|X|_C\otimes|Y|_C,\\(X\Box  Y)\otimes_{\Delta\times\Delta}|\Delta_\EE\Box \Delta_\EE|_C&\cong|X\otimes Y|_C.\end{align*}\end{lma}
\begin{proof}The first isomorphism follows from the facts that the category of simplices of $\Delta^n$ has $id_{[n]}$ as terminal object, and that $|\Delta_\EE^n|_C$ may be identified with the colimit of the functor which takes each simplex $[k]\to[n]$ of $\Delta^n$ to $C^k$. The second isomorphism follows from the first by separating the variables on the left hand side, and using that the tensor of $\EE$ preserves colimits in both variables. The third isomorphism is more involved and is proved in three steps:

(1) Let $d:\AA\to\BB$ be a functor of small categories. Then $d$ induces evident restriction functors $d^*:\EE^\BB\to\EE^\AA$ and $d^*:\EE^{\BB^\op}\to\EE^{\AA^\op}$, the left adjoints of which we denote by $d_!$. Let $B:\BB^\op\to\EE$ and $C:\AA\to\EE$ be functors. Then, $$B\otimes_\BB d_!(C)\cong d^*(B)\otimes_\AA C.$$ Indeed, for each object $E$ of $\EE$, there are bijective correspondences between maps:$$\begin{array}{cccl}B\otimes_\BB d_!(C)&\lra&E&\text{in }\EE\\\hline d_!(C)&\lra&\iHom_\EE(B(-),E)&\text{in }\EE^{\BB}\\\hline C&\lra&d^*\iHom_\EE(B(-),E)&\text{in }\EE^\AA\\\hline C&\lra&\iHom_\EE((d^*B)(-),E)&\text{in }\EE^\AA\\\hline d^*(B)\otimes_{\AA}C&\lra&E&\text{in }\EE\end{array}$$

(2) Specialise to $d:\Delta\to\Delta\times\Delta$, i.e. $B$ is a bisimplicial object of $\EE$, and $C$ is a cosimplicial object of $\EE$. Then, $$d_!(C)\cong d_!(|\Delta_\EE|_C)\cong|d_!(\Delta_\EE)|_C\cong|\Delta_\EE\Box\Delta_\EE|_C.$$Indeed, the left adjoint $d_!$ is actually a left Kan extension which is computed pointwise; in particular, post-composition with a colimit-preserving functor commutes with $d_!$, which yields the second isomorphism above; for the third isomorphism, use that $\Delta_\EE$ is the composite functor $(-)_\EE\circ\Delta_\Sets$, that  $d_!\Delta_\Sets=\Delta_\Sets\Box\Delta_\Sets$, and that $(-)_\EE:\Sets^\Dop\to\EE^\Dop$ is a strong symmetric monoidal functor, whence $(\Delta_\Sets\Box\Delta_\Sets)_\EE=\Delta_\EE\Box\Delta_\EE$.\vspace{1ex}

(3) Specialise to $B=X\Box Y$ and use that $d^*(X\Box Y)=X\otimes Y$.\end{proof}

\begin{prp}\label{EZ}For a cosimplicial object $C$ of a closed symmetric monoidal category $\EE$ with unit $I\cong C^0$, the following three properties are equivalent:
\begin{enumerate}\item $C$ is (strong) monoidal;
\item $|(-)_\EE|_C:\Sets^\Dop\to\EE$ is (strong) symmetric monoidal;
\item There exists a system of Eilenberg-Zilber (iso)morphisms$$EZ^{m,n}_\EE:|\Delta^m_\EE|_C\otimes|\Delta^n_\EE|_C\to|\Delta^m_\EE\otimes\Delta^n_\EE|_C$$ subject to natural associativity, unit and symmetry conditions.\end{enumerate}\end{prp}

\begin{proof}By definitions, (1) implies (2), and (2) implies (3). Condition (3) and Lemma \ref{bi} allow us to endow the realisation functor $|\!-\!|_C$ with a symmetric monoidal structure. The associativity, unit and symmetry conditions of the Eilenberg-Zilber maps translate into the corresponding conditions of the symmetric monoidal structure of $|\!-\!|_C$. Therefore, (3) implies (1).\end{proof}

\begin{rmk}The symmetry and unit conditions of an ``Eilenberg-Zilber system'' simply express compatibility of $EZ^{m,n}_\EE$ with the symmetry and unit of $\EE$. The associativity condition is more subtle to write down explicitly, since it already uses Lemma \ref{bi}. Indeed, the associativity condition for an Eilenberg-Zilber system requires that for any triple $([m],[n],[p])$ of finite ordinals, the following two compositions are equal:

\begin{diagram}[small,silent]&&|\Delta^m_\EE\otimes\Delta^n_\EE|_C\otimes|\Delta^p_\EE|_C&&\\&\ruTo^{EZ^{m,n}_\EE\otimes id_{|\Delta^p_\EE|_C}}&&\rdTo^{(\ref{bi})}&\\|\Delta^m_\EE|_C\otimes|\Delta^n_\EE|_C\otimes|\Delta_p^\EE|_C&&&&|\Delta^m_\EE\otimes\Delta^n_\EE\otimes\Delta^p_\EE|_C\\&\rdTo_{id_{|\Delta^m_\EE|_C}\otimes EZ^{n,p}_\EE}&&\ruTo_{(\ref{bi})}&\\&&|\Delta^m_\EE|_C\otimes|\Delta^n_\EE\otimes\Delta^p_\EE|_C&&\end{diagram}\end{rmk}

\begin{cor}\label{mon}Any unit- and colimit-preserving (strong) symmetric monoidal functor preserves (strong) monoidal cosimplicial objects.\end{cor}

\begin{proof}Consider a functor $\Phi:\EE\to\FF$ and a (strong) monoidal system $C$ in $\EE$. If $\phi$ preserves coproducts, then $(-)_\FF=\Phi\circ(-)_\EE:\Sets^\Dop\to\FF^\Dop$, and this a unit-preserving (strong) symmetric monoidal functor whenever $\Phi$ is. Moreover, if $\phi$ preserves arbitrary colimits then $\phi$ commutes with geometric realisation, in the sense that $|\!-\!|_{\Phi\circ C}$ is naturally isomorphic to the composite functor $\Phi\circ|\!-\!|_C$ . Therefore, condition (2) of Proposition \ref{EZ} gives the conclusion.\end{proof}

\begin{dfn}A cosimplicial object $C$ of a monoidal model category $\EE$ is called \emph{h-monoidal} (resp. \emph{h-constant}) if the realisation functor $|\!-\!|_C:\EE^\Dop\to\EE$ is h-monoidal (resp. if the simplicial operators act as weak equivalences), cf. \ref{functor}. 

A \emph{standard system of simplices} is a cosimplicial object which is h-monoidal, Reedy-cofibrant and h-constant.\end{dfn}

The definition of an h-monoidal cosimplicial object uses the Reedy model structure on $\,\EE^\Dop$, while the definition of a Reedy-cofibrant cosimplicial object uses the Reedy model structure on $\,\EE^\Delta$. We refer to the books of Hirschhorn \cite{Hir} and Hovey \cite{Hov} for the basic properties of Reedy model structures. For our purposes, the following facts suffice: for each Reedy category $\RR$, the category $\EE^\RR$ of functors $\RR\to\EE$, carries a model structure for which the weak equivalences are the natural transformations which are \emph{objectwise weak equivalences}. For each object $X$ of $\EE^\RR$ and each object $r$ of $\RR$, there is a so-called \emph{$r$-th latching map} $L_rX\to X(r)$ in $\EE$; a morphism $X\to Y$ in $\EE^\RR$ is a \emph{Reedy-cofibration} if for all objects $r$ of $\RR$, the induced map $X(r)\cup_{L_rX}L_rY\to Y(r)$ is a cofibration in $\EE$. \emph{Reedy-fibrations} are defined dually by means of so-called \emph{matching maps} $X(r)\to M_rX$. If the cofibrations of $\EE$ are precisely the monomorphisms, the same is true for the Reedy-cofibrations in $\EE^\RR$. For each Reedy category $\RR$, the opposite category $\RR^\op$ is also a Reedy category; the product of two Reedy categories is again a Reedy category in a canonical way. 

The category $\Delta$ is one of the most prominent examples of a Reedy category. For a cosimplicial object $C$, the $[n]$-th latching map is usually denoted by $\partial C^n\to C^n$. For $n=1$, the latching map $\partial C^1\to C^1$ may be identified with the map $C^0\sqcup C^0\to  C^1$ induced by the two coface operators $C^0\dto C^1$. The codegeneracy operator $C^1\to C^0$ is a retraction of these coface operators so that we get a factorisation $C^0\sqcup C^0\to C^1\to C^0$ of the codiagonal like in Definition \ref{interval} of a segment. This leads to the following elementary, but important lemma:

\begin{lma}\label{interval2}The $1$-truncation of a monoidal cosimplicial object is a segment. This segment is cofibrant (resp. an interval) if the cosimplicial object is Reedy-cofibrant (resp. a standard system of simplices).\end{lma}
\begin{proof}The $1$-truncation of the Yoneda-embedding $\Delta\to\Sets^\Dop$ yields a segment $\Delta^0\sqcup\Delta^0\to\Delta^1\to\Delta^0$ in simplicial sets, with binary operation induced by the maximum operation. The Yoneda-embedding is sent to a given system of simplices $C$ in $\EE$ by composition with the functor $|(-)_\EE|_C:\Sets^\Dop\to\EE$, cf. Lemma \ref{bi}. By Proposition \ref{EZ}, the functor $|(-)_\EE|_C$ is unit-preserving and symmetric monoidal. Since any unit-preserving monoidal functor sends segments to segments, the image $C^0\sqcup C^0\to C^1\to C^0$ of $\Delta^0\sqcup\Delta^0\to\Delta^1\to\Delta^0$ is a segment. It follows from the definitions that this segment is cofibrant (resp. an interval) if $C$ is Reedy-cofibrant (resp. a standard system of simplices).\end{proof}

We shall use the following three technical lemmas concerning Reedy model structures on functor categories:

\begin{lma}\label{Reedy}Let $\RR$ be a Reedy category. For each Reedy-cofibrant object $C$ of $\EE^\RR$, the adjunction $(-)\otimes_\RR C:\EE^{\RR^{op}}\leftrightarrows\EE:(-)^{C}$ is a Quillen adjunction.\end{lma}

\begin{proof}We have to show that for any (trivial) fibration $W\to Z$ in $\EE$, the induced map $Z^C\to W^C$ in $\EE^\Dop$ is a (trivial) Reedy-fibration, i.e. for each object $r$ of $\RR$, the map $Z^{C(r)}\to W^{C(r)}\times_{W^{L_rC }}Z^{L_rC}$ should have the right lifting property with respect to trival cofibrations (resp. cofibrations) $A\ito B$ in $\EE$. By exponential transpose, this lifting property translates into the right lifting property of the given (trivial) fibration $W\to Z$ with respect to $A\otimes C(r)\cup_{A\otimes L_rC}B\otimes L_rC\to B\otimes C(r)$; this in turn follows from the pushout-product axiom, since $C$ is Reedy-cofibrant, i.e. $L_rC\to C(r)$ is a cofibration for all objects $r$ of $\RR$.\end{proof}

\begin{lma}\label{Reedy2}Let $\RR,\SS$ be Reedy categories. For Reedy-cofibrant objects $C,D$ of $\EE^\RR,\EE^\SS$, the external tensor-product $C\Box D$ is a Reedy-cofibrant object of $\EE^{\RR\times\SS}$.\end{lma}

\begin{proof}By hypothesis, the latching maps $L_rC\to C(r)$ and $L_sD\to D(s)$ are cofibrations for any object $(r,s)$ of $\RR\times\SS$. Therefore, by the pushout-product axiom, the induced map $L_rC\otimes D(s)\cup_{L_rC\otimes L_sD}C(r)\otimes L_sD\to C(r)\otimes D(s)$ is a cofibration too; the latter may be identified with the latching map $L_{(r,s)}(C\Box D)\to (C\Box D)(r,s)$.\end{proof}

\begin{lma}\label{Reedy3}Let $\Phi:\EE\to\FF$ be a colimit- and cofibration-preserving functor of monoidal model categories. Then, the induced functor $\Phi^\RR:\EE^\RR\to\FF^\RR$ preserves Reedy-cofibrations. In particular, if $\Phi$ is a left Quillen functor, then so is $\Phi^\RR$.\end{lma}

\begin{proof}A (trivial) Reedy-cofibration $X\to Y$ in $\EE^\RR$ is characterised by the property that for all objects $r$ of $\RR$, the induced map $X(r)\cup_{L_rX}L_rY\to Y(r)$ is a (trivial) cofibration. A colimit-preserving functor $\Phi$ takes the latter to $(\Phi X)(r)\cup_{L_r(\Phi X)}L_r(\Phi Y)\to(\Phi Y)(r)$. Therefore, if $\Phi$ preserves (trivial) cofibrations, then $\Phi^\RR$ preserves (trivial) Reedy-cofibrations.\end{proof}

\begin{prp}\label{hEZ}For a cosimplicial object $C$ of a monoidal model category $\EE$, the following three properties are equivalent:\begin{enumerate}\item $C$ is h-monoidal and Reedy-cofibrant;\item $|(-)_\EE|_C:\Sets^\Dop\to\EE$ is h-monoidal and preserves cofibrations;\item $C$ is Reedy-cofibrant and comes equipped with a system of Eilenberg-Zilber morphisms $EZ^{m,n}_\EE$ (like in Proposition \ref{EZ}) which are weak equivalences.\end{enumerate}\end{prp}

\begin{proof}If we assume (1), Lemma \ref{Reedy} implies that the realisation functor $|\!-\!|_C$ is an h-monoidal left Quillen functor; the composite functor $|(-)_\EE|_C$ is therefore cofibration-preserving and h-monoidal, whence (1) implies (2). Domain and codomain of the Eilenberg-Zilber maps are cofibrant objects of $\EE$, since they belong to the image of the functor $|(-)_\EE|_C$, whence (2) implies (3). In order to prove (3) implies (1), it remains to be shown that, for cofibrant objects $X$ and $Y$ of $\EE^\Dop$, the structural map $|X|_C\otimes|Y|_C\to|X\otimes Y|_C$, given by Proposition \ref{EZ}, is a weak equivalence. By Lemma \ref{bi}, this map may be identified with a map of coends\begin{diagram}(X\Box Y)\otimes_{\Delta\times\Delta}(|\Delta_\EE|_C\Box|\Delta_\EE|_C)&\rTo^{id_{X\Box Y}\otimes_{\Delta\times\Delta} EZ^{-,-}_\EE}&(X\Box Y)\otimes_{\Delta\times\Delta}|\Delta_\EE\Box \Delta_\EE|_C.\end{diagram}The right hand side components are precisely the Eilenberg-Zilber maps which are supposed to be weak equivalences. By Lemma \ref{Reedy}, it suffices therefore to show that $X\Box Y$ is Reedy-cofibrant in $\EE^{(\Delta\times\Delta)^{op}}$, and that $|\Delta_\EE|_C\Box |\Delta_\EE|_C$ and $|\Delta_\EE\Box \Delta_\EE|_C$ are Reedy-cofibrant in $\EE^{\Delta\times\Delta}$. For the first two objects, this follows from from Lemmas \ref{Reedy2} and \ref{bi}. 

For the third object, observe first that $\Delta_\EE$ is the image of the Yoneda-embedding $\Delta_\Sets$ under the functor $(-)_\EE$; Lemma \ref{Reedy3} thus implies that $\Delta_\EE$ is Reedy-cofibrant in $(\EE^\Dop)^\Delta$. Therefore, by Lemma \ref{Reedy2}, $\Delta_\EE\Box\Delta_\EE$ is Reedy-cofibrant in $(\EE^\Dop)^{\Delta\times\Delta}$. Finally, the image of the latter under the cofibration-preserving functor $|\!-\!|_C$ is precisely the object $|\Delta_\EE\Box \Delta_\EE|_C$ under consideration, which is thus Reedy-cofibrant in $\EE^{\Delta\times\Delta}$ by Lemma \ref{Reedy3}.\end{proof}

\begin{cor}\label{hmon}Any colimit- and cofibration-preserving h-monoidal functor preserves h-monoidal Reedy-cofibrant cosimplicial objects.\end{cor}
\begin{proof}This follows from Corollary \ref{mon}, Lemma \ref{Reedy3} and Proposition \ref{hEZ}.\end{proof}

\begin{prp}\label{left}For a cosimplicial object $C$ of a monoidal model category $\EE$, the following three properties are equivalent:\begin{enumerate}\item $C$ is a standard system of simplices;\item $|(-)_\EE|_C:\Sets^\Dop\to\EE$ is an h-monoidal left Quillen functor;\item $C$ is h-monoidal, Reedy-cofibrant, and $C^1\to C^0$ is a weak equivalence.\end{enumerate}\end{prp}

\begin{proof}That (1) implies (3) is immediate. (2) implies (1), since Proposition \ref{hEZ} shows that under hypothesis (2), $C$ is h-monoidal and Reedy-cofibrant; moreover, as the operators of $\Delta$ act as weak equivalences between the representable simplices in $\Sets^\Dop$ and any left Quillen functor preserves weak equivalences between cofibrant objects, $C$ is h-constant. 

It remains to prove (3) implies (2). Note that the functor $|(-)_\EE|_C$ has a right adjoint mapping the object $X$ to the simplicial set $\Sing_C(X)$ given by $\Sing_C(X)_n=\EE(C^n,X)$. Thus, by Proposition \ref{hEZ}, we only need to show that the functor $|(-)_\EE|_C$ preserves trivial cofibrations. For this, it suffices to consider a complete set of generating trivial cofibrations for simplicial sets; such a set is given by the inclusions $\Delta^n\times\Delta^0\cup_{\partial\Delta^n\times\Delta^0}\partial\Delta^n\times\Delta^1\hookrightarrow\Delta^n\times\Delta^1$ for $0,1:\Delta^0\dto\Delta^1$ and $n\geq 0.$ It follows from Reedy's patching lemma that the h-monoidal functor $|(-)_\EE|_C$ sends these inclusions to morphisms weakly equivalent to $C^n\otimes C^0\cup_{\partial C^n\otimes C^0}\partial C^n\otimes C^1\to C^n\otimes C^1$. The latter are trivial cofibrations by an application of the pushout-product axiom, since $\partial C^n\to C^n$ is a cofibration, and $0,1:C^0\dto C^1$ are trivial cofibrations.\end{proof}

\begin{cor}\label{hfc}Any h-monoidal left Quillen functor preserves standard systems of simplices.\end{cor}



\begin{rmk}\label{framing}Any standard system of simplices $\,C$ in $\EE$ induces so-called \emph{cosimplicial} resp. \emph{simplicial framings} $-\otimes C$ resp. $(-)^{C}$ for all objects of $\EE$, cf. \cite[Chapter 5]{Hov}, \cite[Chapter 17]{Hir}. This means that we can define simplicial mapping spaces $$\EE(X\otimes C, Y)\cong\EE(X,Y^C)\cong\EE(C,Y^X),$$ which have the ``correct homotopy type'' whenever $X$ is cofibrant and $Y$ is fibrant. Moreover, condition (2) of Proposition \ref{left} implies that Quillen's axiom SM7 holds for these mapping spaces. Observe however that we need further structure in order to get an enrichment over simplicial sets. Indeed, it is sufficient to require that the standard system of simplices be also comonoidal in a sense analogous to Definition \ref{defmon} (no symmetry constraint is needed); the composition of simplicial mapping spaces is then given by $C\to C\otimes C\to Z^Y\otimes Y^X\to Z^X$. Of course, this condition is automatically fulfilled if the tensor is the categorical product. Even more structure is needed to get a \emph{simplicial model structure} on $\EE$, since the latter requires $\EE$ to be tensored and cotensored over the category of simplicial sets; this is precisely the case when the standard system of simplices is strong monoidal (and not only $h$-monoidal).\end{rmk}

\begin{exms}\label{standard}The standard system of simplices for simplicial sets is given by the Yoneda-embedding $\Delta\to\Sets^\Dop$. This system is strong symmetric monoidal.

The standard system of simplices for compactly generated spaces is obtained as the geometric realisation of the standard system for simplicial sets. The geometric realisation functor is a strong symmetric monoidal left Quillen functor, therefore Proposition \ref{left} applies. Similarly, the suspension spectrum functor is a strong symmetric monoidal left Quillen functor from simplicial sets to symmetric spectra (or any other symmetric monoidal model for stable homotopy, cf. \cite{MMSS}), whence a strong monoidal standard system of simplices for symmetric spectra. 

The normalised $R$-chain functor $\NN_*^R$ from simplicial sets to $R$-chain complexes is an h-monoidal left Quillen functor, whence an h-monoidal standard system of simplices for $R$-chain complexes. We borrowed the terminology of an Eilenberg-Zilber system (cf. Propositions \ref{EZ} and \ref{hEZ}) from this special case, since the corresponding ``shuffle'' maps $\NN_*^R(\Delta^m)\otimes_R\NN_*^R(\Delta^n)\to\NN_*^R(\Delta^m\times\Delta^n)$ have been introduced by Eilenberg and Zilber. There are canonical retractions of these shuffle maps introduced by Alexander and Whitney. The latter provide a comonoidal structure for $\NN_*^R$. Therefore, $\Ch(R)$ is actually enriched over simplicial sets without being an honest simplicial model category, cf. Remark \ref{framing}.\end{exms}

\vspace{5ex}

\noindent{\sc Universit\'e de Nice, Laboratoire J.-A. Dieudonn\'e, Parc Valrose, 06108 Nice Cedex, France.}\hspace{2em}\emph{E-mail:} cberger$@$math.unice.fr\vspace{2ex}

\noindent{\sc Mathematisch Instituut, Postbus 80.010, 3508 TA Utrecht, The Ne-therlands.}\hspace{2em}\emph{E-mail:} moerdijk$@$math.uu.nl

\end{document}